\documentclass{jams-l}
\usepackage{geompsfi}

\title[The Weil-Petersson metric and volumes of convex cores]{The Weil-Petersson metric 
and volumes of 3-dimensional hyperbolic convex cores}
\date{October 30, 2001}
\author{Jeffrey F. Brock}
\address{\scriptsize Mathematics Department, University of Chicago, 5734 S. University
Ave., Chicago, IL 60637} 
\email{brock@math.uchicago.edu}
\thanks{Research partially supported by NSF
grant DMS-0072133 and an NSF postdoctoral fellowship.}

\newtheorem{theorem}{Theorem}[section]

\renewcommand{\bold}[1]{\medskip \noindent {\bf \boldmath #1
                        }\nopagebreak[4]}

\renewcommand{\qed}{\nopagebreak[4]\begin{flushright} 
\rule{2mm}{2.5mm} \end{flushright}\pagebreak[2]}

\newcommand{\cx}{{\mathbb C}}
\newcommand{\half}{{\mathbb H}}
\newcommand{\integers}{{\mathbb Z}}
\newcommand{\natls}{{\mathbb N}}

\newcommand{\reals}{{\mathbb R}}

\newcommand{\makefig}[3]{
	\begin{figure}[htbp]
        \refstepcounter{figure}
	\label{#2}
        \begin{center}~
		#3~\\
		\medskip
                {\sf Figure \thefigure.  #1}
        \end{center}
	\medskip
	\end{figure}
}

\newenvironment{pf*}[1]{%
 \begin{proof}[#1]%
}{ 
 \end{proof}
}

\newcommand{\ital}[1]{\medskip \noindent {\em #1 }\nopagebreak[4]}

\newcommand{\bdry}{\partial}

\newcommand{\closure}{\overline}
\newcommand{\compos}{\circ}

\newcommand{\disjunion}{\sqcup}

\usepackage{amssymb}
\newcommand{\nullset}{\varnothing}

\newcommand{\st}{\; | \;}         

\newcommand{\wt}{\widetilde}
\newcommand{\wh}{\widehat}

\newcommand{\chat}{\widehat{\cx}}

\newcommand{\zed}{\integers}

\newcommand{\area}{\mbox{\rm area}}

\newcommand{\core}{\mbox{\rm core}}

\newcommand{\diam}{\mbox{\rm diam}}

\newcommand{\interior}{\mbox{\rm int}}

\newcommand{\Isom}{\mbox{\rm Isom}}

\newcommand{\mass}{\mbox{\rm mass}}

\newcommand{\Mod}{\mbox{\rm Mod}}

\newcommand{\PSL}{\mbox{\rm PSL}}

\newcommand{\sh}{\calS \calH}
\newcommand{\Sing}{\mbox{\rm Sing}}

\newcommand{\Teich}{\mbox{\rm Teich}}

\newcommand{\vol}{\mbox{\rm vol}}

\newtheorem{prop}[theorem]{Proposition}
\newtheorem{lem}[theorem]{Lemma}

\newtheorem{defn}[theorem]{Definition}

\newcommand{\calA}{{\mathcal A}}
\newcommand{\calB}{{\mathcal B}}
\newcommand{\calC}{{\mathcal C}}

\newcommand{\calH}{{\mathcal H}}
\newcommand{\calI}{{\mathcal I}}

\newcommand{\calM}{{\mathcal M}}
\newcommand{\calN}{{\mathcal N}}

\newcommand{\calP}{{\mathcal P}}

\newcommand{\calS}{{\mathcal S}}
\newcommand{\calT}{{\mathcal T}}

\newcommand{\calV}{{\mathcal V}}

\usepackage{euscript}

\newcommand{\eG}{{\EuScript G}}

\newcommand{\eS}{{\EuScript S}}
\newcommand{\eT}{{\EuScript T}}

\newcommand{\eV}{{\EuScript V}}

\begin{document}

\begin{abstract}
We present a coarse interpretation of the Weil-Petersson distance
$d_{\rm WP}(X,Y)$ between two finite area hyperbolic Riemann surfaces
$X$ and $Y$ using a graph of pants decompositions introduced by
Hatcher and Thurston.  The combinatorics of the pants graph reveal a
connection between Riemann surfaces and hyperbolic 3-manifolds
conjectured by Thurston: the volume of the convex core of the
quasi-Fuchsian manifold $Q(X,Y)$ with $X$ and $Y$ in its conformal
boundary is comparable to the Weil-Petersson distance $d_{\rm
WP}(X,Y)$.  In applications, we relate the Weil-Petersson
distance to the Hausdorff dimension of the limit set and
the lowest eigenvalue of the Laplacian for $Q(X,Y)$, and give a new
finiteness criterion for geometric limits.
\end{abstract}

\maketitle

\makefig{The lift to $\half^3$ of a quasi-Fuchsian convex
core boundary component.
}{figure:pleated:plane}{
\begin{center}~
\psfig{file=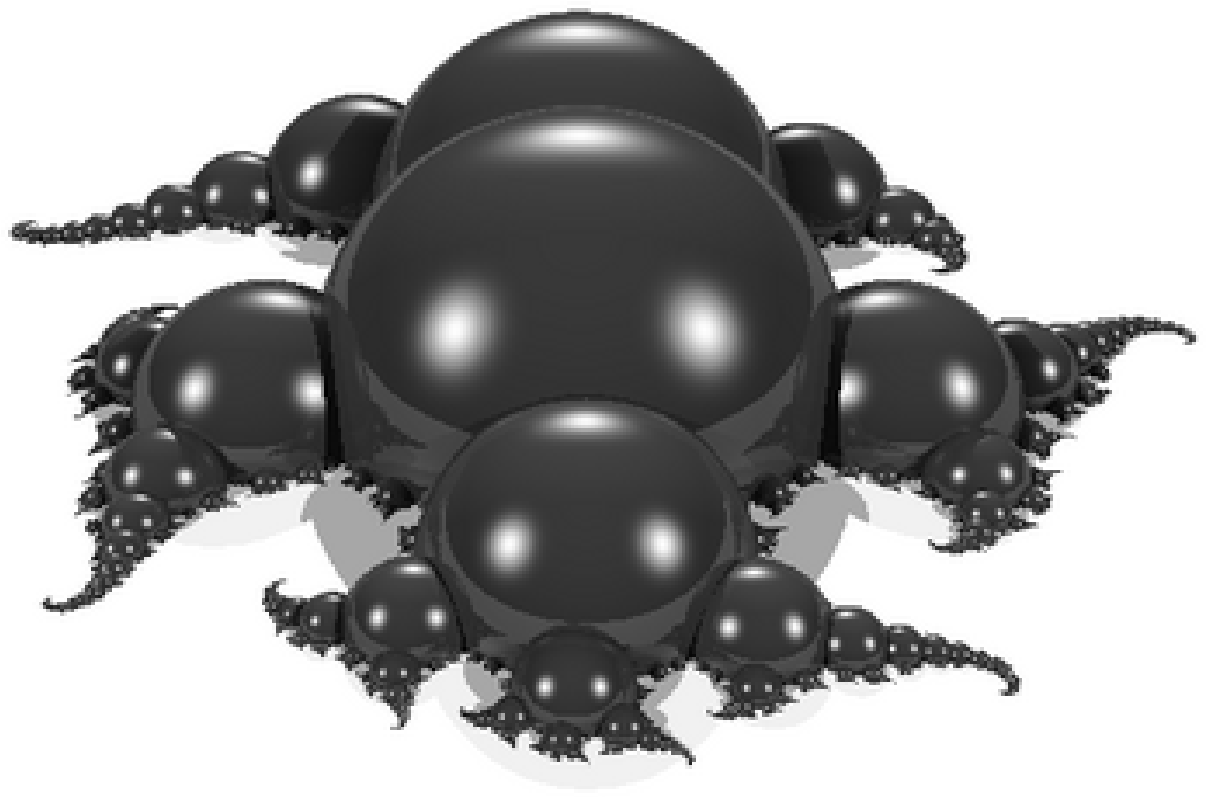,height=2.5in}
\end{center}
}

\section{Introduction}
\label{section:introduction}  
Recent insights into the combinatorial geometry of
Teichm\"uller space have shed new light on fundamental questions in
hyperbolic geometry in 2 and 3 dimensions. Paradoxically, a coarse
perspective on Teichm\"uller space appears to refine the analogy of
Teichm\"uller geometry with the internal geometry of hyperbolic
3-manifolds first introduced and pursued by W. Thurston.

In this paper we develop such a coarse perspective on the Weil-Petersson
metric on Teichm\"uller space by relating it to a graph of
pair-of-pants decompositions of surfaces introduced by Hatcher and
Thurston.  This viewpoint generates a 
new connection between the Weil-Petersson geometry of Teichm\"uller
space and the geometry of the convex core of a hyperbolic 3-manifold.

\smallskip
For simplicity, let $S$ be a closed oriented surface of negative Euler
characteristic.  A {\em pants decomposition} of $S$ is a maximal
collection of distinct isotopy classes of pairwise disjoint essential
simple closed curves on $S$.  We say two distinct pants decompositions
$P$ and $P'$ are related by an {\em elementary move} if $P'$ can be
obtained from $P$ by replacing a curve $\alpha \in P$ by a curve
$\beta$ intersecting $\alpha$ minimally (see
Figure~\ref{figure:move}).

One obtains the {\em pants graph} ${\bf P}(S)$ by making each pants
decomposition a vertex and joining two pants decompositions differing
by an elementary move by an edge.  
Setting the length of each edge to
1, ${\bf P}(S)$ becomes a metric space.  
We find the graph ${\bf
P}(S)$ provides a combinatorial model for the coarse geometry of the
Weil-Petersson metric:
\begin{theorem}  
The graph ${\bf P}(S)$ is naturally quasi-isometric to Teichm\"uller space
with the Weil-Petersson metric.
\label{theorem:combinatorics}
\end{theorem}

The connection to hyperbolic 3-manifolds is simple to describe.  By a
theorem of Bers, a pair of points $(X,Y) \in \Teich(S) \times
\Teich(S)$ naturally determines a quasi-Fuchsian hyperbolic 3-manifold
$Q(X,Y) \cong S \times \reals$ with $X$ and $Y$ in its conformal
boundary at infinity.  Its {\em convex core}, denoted $\core(Q(X,Y))$,
is the smallest convex subset of $Q(X,Y)$ carrying its fundamental
group.  The convex core is itself homeomorphic to $S \times I$ and
carries all the essential geometric information about the manifold $Q(X,Y)$.

Because $Q(X,Y)$ is obtained from the pair $(X,Y)$ by an
analytic process (Bers's {\em simultaneous uniformization}), it is a
central challenge in the study of hyperbolic  
3-manifolds to understand the geometry of $Q(X,Y)$ purely 
in terms of the geometry of $X$ and $Y$.
Our main theorem proves a conjecture of Thurston that the following
fundamental connection exists between 
convex core
{\em volume} and the Weil-Petersson distance.
\begin{theorem}
The volume of the convex core of $Q(X,Y)$ is comparable to the
Weil-Petersson distance $d_{\rm WP}(X,Y)$.
\label{theorem:main}
\end{theorem}

Here, {\em comparability} means that two quantities are equal up to  
uniform additive and multiplicative error: i.e.
there are constants
$K_1 >1$ and $K_2 >0$ depending only on $S$ so that for any $(X,Y) \in 
\Teich(S) \times \Teich(S)$ we have
$$\frac{d_{\rm WP}(X,Y)}{K_1} - K_2 \le 
\vol(\core(Q(X,Y))) \le
K_1 d_{\rm WP}(X,Y) + K_2.$$ Throughout the paper we will use the
contraction $\vol(X,Y) = \vol(\core(Q(X,Y)))$ and the notation
$\asymp$ to denote the comparability of two quantities; then
Theorem~\ref{theorem:main} becomes $$d_{\rm WP}(X,Y) \asymp
\vol(X,Y).$$

The volume of the convex core of a complete hyperbolic 3-manifold $M =
\half^3 /\Gamma$ is
directly related to the lowest eigenvalue of the Laplacian on $M$ as
well as the Hausdorff dimension of the {\em limit set} $\Lambda(\Gamma)
\subset \chat$, namely the complement of the invariant {\em domain of
discontinuity} $\Omega(\Gamma) \subset \chat$ where the action of
Kleinian covering group $\Gamma \subset
\PSL_2(\cx)$ for $M$ is properly discontinuous (see
Figure~\ref{fig:limitsets} for two examples of limit sets\footnote{We have employed computer programs of Curt McMullen in our generation of
Figures~\ref{figure:pleated:plane} and~\ref{fig:limitsets}.}).

\makefig{Limit sets for two quasi-Fuchsian groups.
}
{fig:limitsets}{
\begin{center}~
\psfig{file=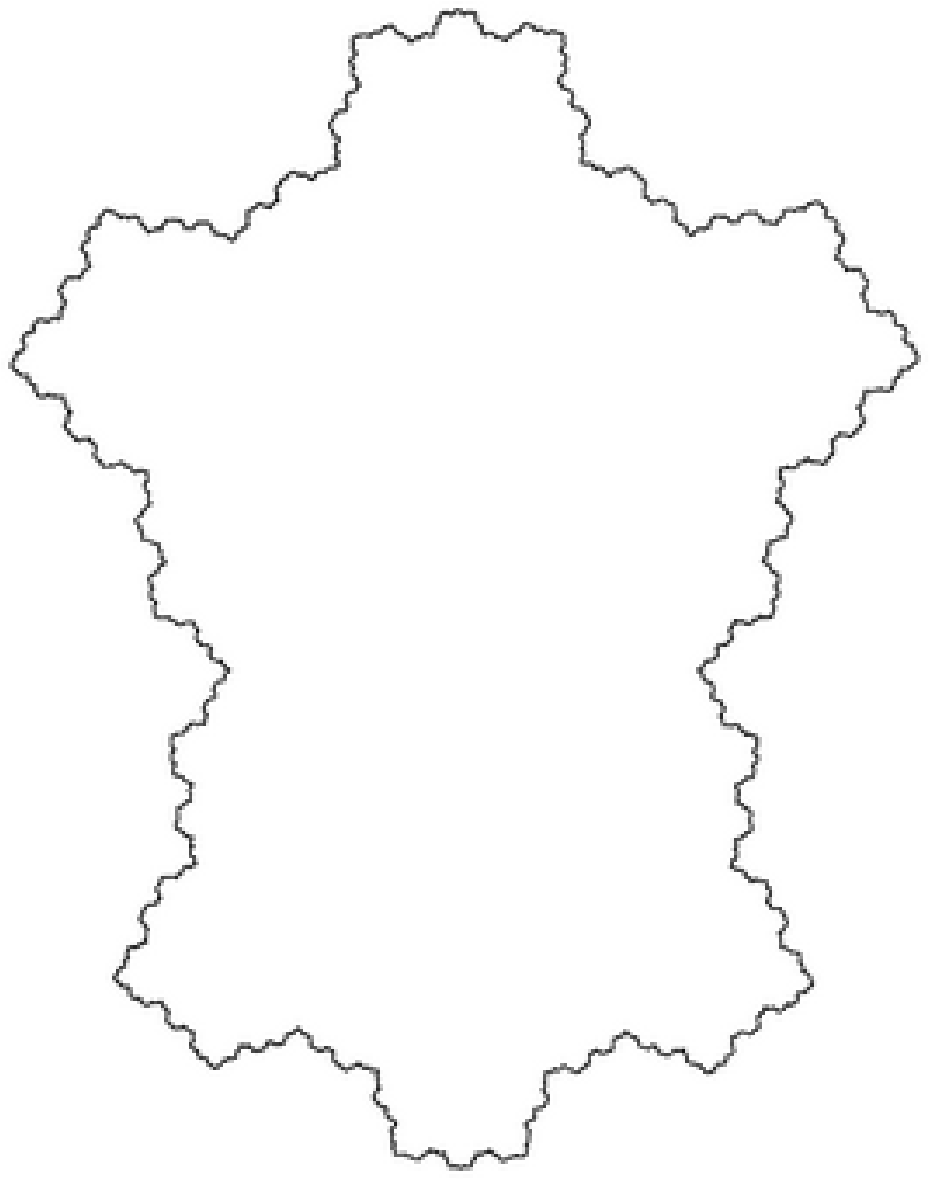,height=2.5in}~
\psfig{file=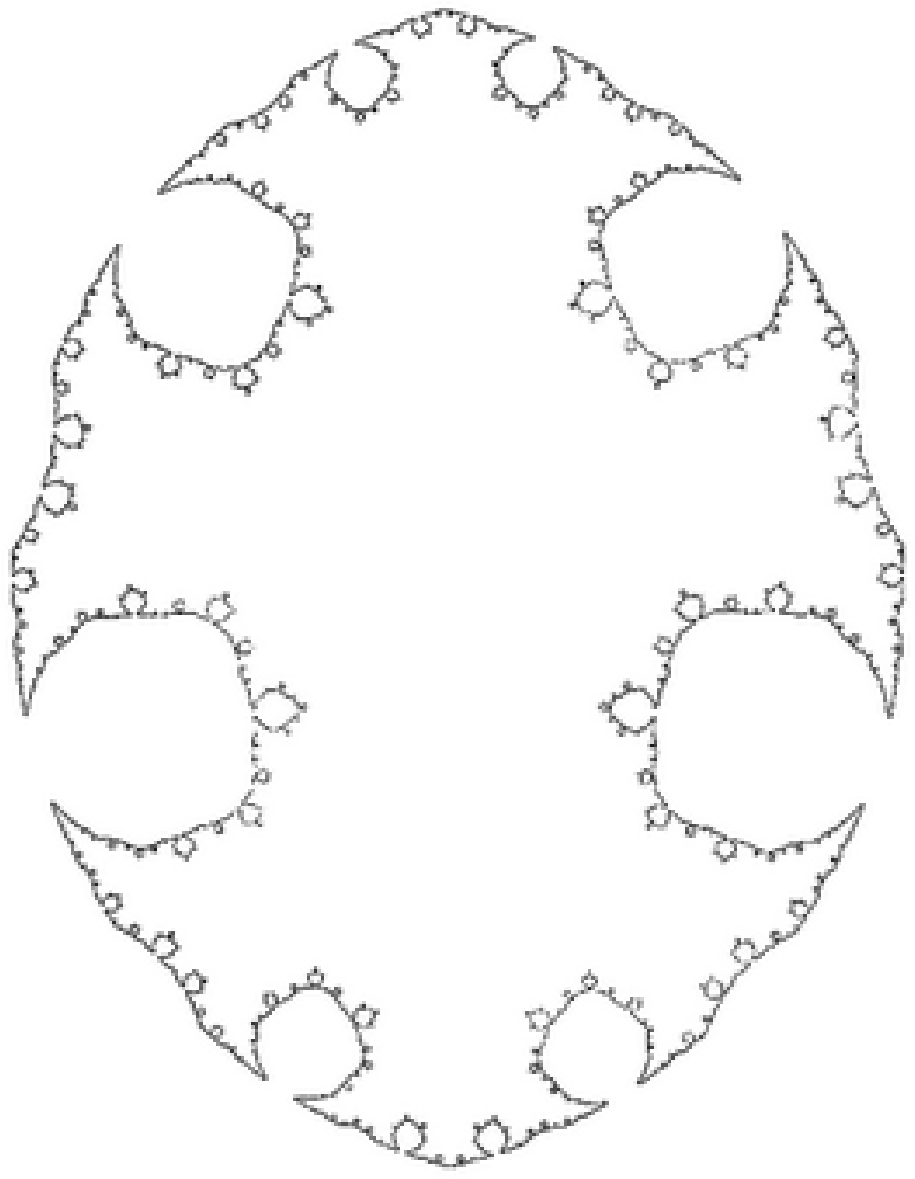,height=2.5in}
\end{center}
}

As an immediate application, Theorem~\ref{theorem:main} implies the
following new relationship between these analytic invariants and the
Weil-Petersson distance.  Let $\lambda_0(X,Y)$ denote the lowest
eigenvalue of the Laplacian on the quasi-Fuchsian hyperbolic
3-manifold $Q(X,Y) = \half^3/ \Gamma(X,Y)$ and let
$D(X,Y)$ denote the Hausdorff dimension of the limit set of
$\Gamma(X,Y)$.
\begin{theorem}
Given $S$ there are constants $K>0$, $C_1$, $C_2$, $C_3$, and $C_4 >1$ so
that if $d_{\rm WP}(X,Y) > K$ 
then 
$$\frac{C_1}{d_{\rm 
WP}(X,Y)^2} \le \lambda_0(X,Y) \le
\frac{C_2}{d_{\rm WP}(X,Y)},$$
and
$$
2 - \frac{C_3}{d_{\rm WP}(X,Y)}
\le  D(X,Y)
\le 
2 - \frac{C_4}{d_{\rm 
WP}(X,Y)^2} 
.$$
\label{theorem:analytic}
\end{theorem}

\bold{Proof:}  The relation $\lambda_0(X,Y) = D(X,Y)(2 -D(X,Y))$
follows from a general result by D. Sullivan (see
\cite[Thm. 2.17]{Sullivan:lambda0}), after applying Bowen's
Theorem \cite{Bowen:HD} that $D(X,Y) \ge 1$ with equality if and only
if $X = Y$. 

Theorem~\ref{theorem:main} may be rephrased to claim the existence of
$K$, $K'$ so that for $d_{\rm WP}(X,Y) > K$ we have
$$\frac{d_{\rm WP}(X,Y)}{K'} \le \vol(X,Y) \le K' d_{\rm WP}(X,Y).$$
The theorem then follows from 
the double inequality
$$\frac{c_1}{\vol(X,Y)^2} \le \lambda_0(X,Y) \le
\frac{c_2}{\vol(X,Y)}$$
(see \cite[Main Thm.]{Burger:Canary:lambda0}
and \cite[Thm. A]{Canary:lambda0})
after collecting constants.
\qed

\bold{The pants graph.} 
Since Theorem~\ref{theorem:main} relies directly on
Theorem~\ref{theorem:combinatorics} we detail our coarse perspective
on the Weil-Petersson metric.
\makefig{Elementary moves on pants
decompositions.}{figure:move}{\psfig{file=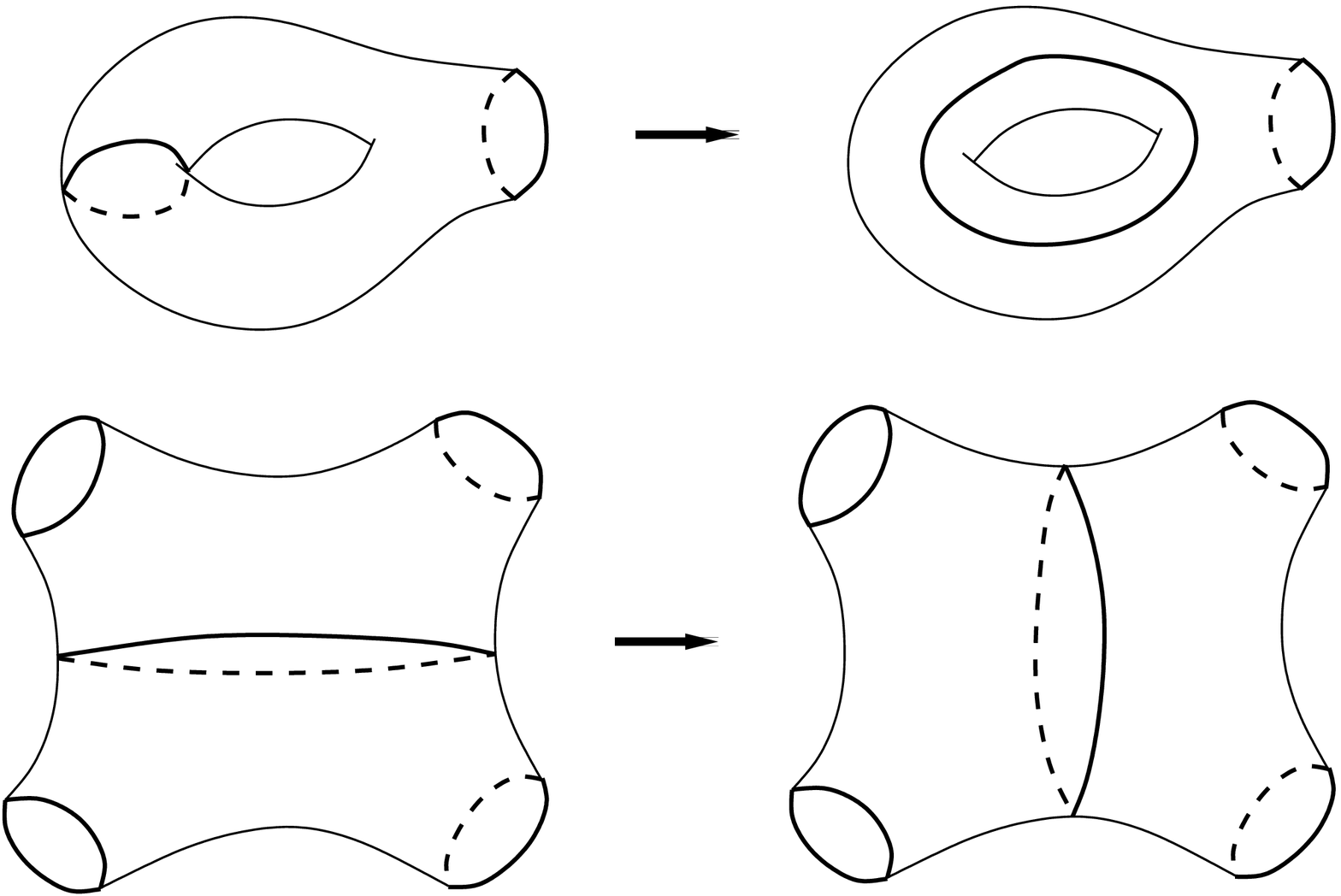,height=2in}} 

To describe the nature of the quasi-isometry between the graph ${\bf
P}(S)$ and the Weil-Petersson metric, we recall that by a theorem of
Bers, there is a constant $L>0$ depending only on $S$ so that for each
$X \in
\Teich(S)$ there is a pants decomposition $P$ so that 
$$\ell_X(\alpha) < L \ \ \ \text{for each} \ \ \ \alpha \in P,$$
where $\ell_X(\alpha)$ denotes the length of the geodesic
representative of $\alpha$ in the hyperbolic metric on $X$.

If $P$ is a pants decomposition,
we denote by $V_\ell(P)$ the sub level set
$$V_\ell(P) = \left\{ X \st \max_{\alpha \in P}\left( \ell_X(\alpha)
\right) < \ell \right\}.$$ 
Then Bers's theorem guarantees that the sub level sets $V(P) = V_L(P)$ cover
Teichm\"uller space.  

Coarsely, the Weil-Petersson distance records the configuration of the
sub level sets $V(P)$ with ${\bf P}(S)$ as its model.
Indeed, 
if $P_X$ and $P_Y$ are pants decompositions for which
$X \in V(P_X)$ and $Y \in V(P_Y)$,
then the quasi-isometry of Theorem~\ref{theorem:combinatorics} arises
from the comparability
$$d_{\bf P}(P_X,P_Y) \asymp d_{\rm WP}(X,Y)$$
where $d_{\bf P}(P_X,P_Y)$ is the minimal number of elementary moves
required to travel from $P_X$ to $P_Y$ in ${\bf P}(S)$.

We briefly outline our proof of
Theorem~\ref{theorem:combinatorics} and its application to
Theorem~\ref{theorem:main}. 

\bold{Outline of the proof of Theorem~\ref{theorem:combinatorics}.}
Let ${\bf P}^0(S)$ denote the vertex set of ${\bf P}(S)$, and let
$Q \colon {\bf P}^0(S) \to \Teich(S)$ be any map so that $Q(P)
\in V(P)$.  Applying work of S. Wolpert and H. Masur there is a uniform
constant $D>0$ so that the Weil-Petersson diameter of $V(P)$ satisfies
$$\diam_{\rm WP}(V(P)) < D$$ for all $P \in {\bf P}^0(S)$.  A simple
argument shows that for any two pants decompositions $P_1$ and $P_2$ differing
by an elementary move there is a single Riemann surface $X$ on 
which all curves in $P_i$ have length less than $L$.
Thus $V(P_1) \cap V(P_2)$ is non-empty, and $Q$ is $2D$-Lipschitz.

Given two pants decompositions $P_1$ and $P_2$ for which $V_{2L}(P_1)
\cap V_{2L}(P_2) \not= \emptyset$, there is an $X \in V_{2L}(P_1) \cap
V_{2L}(P_2)$: i.e. a single Riemann surface on which each curve in
$P_1 \cup P_2$ has length at most $2L$.  It follows that there is a
uniform $C$ so that the total intersection number satisfies
$$i(P_1,P_2) \le C,$$ which in turn provides a uniform bound to the
distance $d_{\bf P}(P_1,P_2)$ in ${\bf P}(S)$.

A compactness argument shows that each $X \in V(P) = V_L(P)$ lies a uniformly
definite distance from $\bdry V_{2L}(P)$.  Thus, a unit-length
Weil-Petersson geodesic can always be covered by a uniform number of
sub level sets $V_{2L}(P)$.  It follows that any pair of pants
decompositions $P$ and $P'$ for which $V(P)$ and $V(P')$ contain the
endpoints of a unit length Weil-Petersson geodesic, $P$ and $P'$ have
uniformly bounded distance in ${\bf P}(S)$, and the theorem follows.

\medskip

\bold{Outline of the proof of Theorem~\ref{theorem:main}.}
The proof has two parts.

\ital{Bounding volume from below:}
The bound below of core volume in terms of the Weil-Petersson distance
begins with an interpolation through the convex core $$h_t \colon Z_t
\to \core(Q(X,Y))$$ of 1-Lipschitz maps of hyperbolic surfaces.
It follows that for each essential simple closed curve $\alpha$ we
have $$\inf_t \ell_{Z_t}(\alpha) \ge
\ell_{Q(X,Y)}(\alpha).$$

The path $Z_t$, then, only passes through sets $V(P)$ for which each
element in $P$ has length less than $L$ in $Q(X,Y)$.  Applying recent
work of Masur and Minsky, we show if a sequence $\{P_1, \ldots,
P_n\}$ of pants decompositions is built from $N$ curves and makes
bounded jumps, i.e.
$$d_{\bf P}(P_j,P_{j+1}) < k,$$
then its ends satisfy the bound $d_{\bf P}(P_1, P_n) < K_0N$, where $K_0$ depends only 
on $k$ and $S$.

The Margulis lemma forces closed geodesics with length less than $L$
in $Q(X,Y)$ that represent different isotopy classes to be uniformly
equidistributed through the convex core.  Since each such representative
makes a definite contribution to core volume, the lower bound follows.

\ital{Bounding volume from above.} Given pants decompositions $P_X$
and $P_Y$ so 
that $X \in V(P_X)$ and $Y \in V(P_Y)$, and a geodesic $G \subset {\bf 
P}(S)$ joining $P_X$ to $P_Y$, we consider the closed geodesics
$${\rm spin}(G) = \{ \alpha^* \st \alpha \in P \ \text{for} \ P \in G
\}$$
where $\alpha^*$ denotes the geodesic representative of $\alpha$ in $Q(X,Y)$.
We build a straight triangulation $\calT$ of all but a uniformly bounded
volume portion of  $\core(Q(X,Y))$ so that vertices of $\calT$ lie on
$\alpha^* \in  {\rm spin}(G)$, the so-called {\em spinning geodesics}.

Our triangulation has the property that all but constant times $d_{\bf
P}(P_X,P_Y)$ of the tetrahedra in $\calT$ have at least one edge in a
spinning geodesic $\alpha^*$.  We
then use a spinning trick: by homotoping the vertices around the
geodesics in ${\rm spin}(G)$ keeping the triangulation straight, all
tetrahedra with an edge in any $\alpha^*$ can be made to have
arbitrarily small volume.  

Since there is an {\em a priori} bound to
the volume of a tetrahedron in $\half^3$, 
the remaining tetrahedra have
uniformly bounded volume. 
The theorem then follows from the
comparability
$d_{\rm WP}(X,Y)  \asymp d_{\bf P}(P_X,P_Y).$

\bold{Geometrically finite hyperbolic 3-manifolds.}  We remark that
simple generalizations of these techniques may be employed to obtain
estimates for core volume of hyperbolic 3-manifolds that are
not quasi-Fuchsian once the appropriate version of Weil-Petersson
distance is defined.  For example, given a hyperbolic 3-manifold
$M_\psi$ that fibers over the circle with monodromy $\psi$, the volume
of $M_\psi$ is comparable to the Weil-Petersson translation distance
of $\psi$ (with constants depending only on the topology of the
fiber).  We take up these generalizations in \cite{Brock:3ms1}.

\bold{Algebraic and geometric limits.}  As an application of
Theorem~\ref{theorem:main}, boundedness of the Weil-Petersson distance
$d_{\rm WP}(X_k,Y_k)$ for sequences predicts the geometric
finiteness of the geometric limit of $Q(X_k,Y_k)$.

The space $QF(S)$ of all quasi-Fuchsian hyperbolic 3-manifolds lies in
the space $AH(S)$ of all complete hyperbolic 3-manifolds $M$ marked by
homotopy equivalences $(h \colon S \to M)$ so that $h_*$ sends
peripheral elements of $\pi_1(S)$ to parabolic elements of $\pi_1(M)$.
The space $AH(S)$ carries the {\em algebraic topology} or the
compact-open topology on the induced representations $h_* \colon
\pi_1(S) \to \Isom^+(\half^3)$ up to conjugacy.

In an algebraically convergent sequence $\{(h_k \colon S \to M_k)\}$
in $AH(S)$, normalizing the induced representations $\rho_k = (h_k)_*$
to converge on generators one may always extract a subsequence so that
the groups $\rho_k(\pi_1(S)) = \Gamma_k$ converge in the
Gromov-Hausdorff topology on discrete subgroups of $\Isom^+(\half^3)$, or
{\em geometrically}, to a limit $\Gamma_G$.  A central issue in the
deformation theory of hyperbolic 3-manifolds is to understand the {\em
geometric limit} $N_G = \half^3 /\Gamma_G$.

Applying Theorem~\ref{theorem:main}, we obtain the following
criterion:
\begin{theorem}
Let $Q(X_k,Y_k) \to Q_\infty$ be an algebraically convergent sequence in
$AH(S)$ with geometric limit $N_G$.  Then $N_G$ is geometrically
finite if and only if there is a $K>0$ for which
$$d_{\rm WP}(X_k,Y_k) < K$$ for all $k$.
\label{theorem:geometric:limit}
\end{theorem}
Note that geometric finiteness of $N_G$ implies geometric finiteness
of $Q_\infty$ but {\em not} conversely.


\bold{History and references.}  The fundamental properties of the
Weil-Petersson metric we use are discussed in
\cite{Wolpert:noncompleteness}, \cite{Wolpert:Nielsen},
\cite{Wolpert:universal} 
and \cite{Masur:WP}.  
The pants graph is the 1-skeleton of the {\em pants
complex}, introduced in \cite{Hatcher:Thurston:pants} (see also
\cite{Hatcher:pants}) which is there proven to be connected.
The relation of the pants graph to the Weil-Petersson metric is
similar in spirit to the relative hyperbolicity theorem for
Teichm\"uller space of \cite{Masur:Minsky:CCI} where the (related)
{\em complex of curves} is shown to be quasi-isometric to the {\em
electric} Teichm\"uller space, and to be Gromov-hyperbolic (the pants
complex and the Weil-Petersson metric are {\em not} in general
Gromov-hyperbolic \cite{Brock:Farb:rank}).  For more on quasi-Fuchsian
manifolds and their algebraic and geometric limits, see
\cite{Thurston:book:GTTM}, \cite{Bers:simunif}, 
\cite{McMullen:book:RTM},  \cite{McMullen:iter}, 
\cite{Brock:iter}, 
and \cite{Otal:book:fibered}.

\bold{Plan of the paper.}
After discussing the fundamental work of S. Wolpert and H. Masur  on
the Weil-Petersson metric that will 
serve as our jumping off point in section~\ref{section:preliminaries}, we 
prove the comparability of Weil-Petersson distance and pants distance
(Theorem~\ref{theorem:combinatorics}) in section~\ref{section:combinatorics}. 
We then establish the lower bound on $\vol(X,Y)$ in terms of
the distance $d_{\bf P}(P_X,P_Y)$ in section~\ref{section:lower}.
Section~\ref{section:upper} applies the
combinatorics of pants decompositions along a geodesic $G 
\subset {\bf P}(S)$ joining $P_X$ to $P_Y$ to bound volume from
above in terms of pants distance.  Theorem~\ref{theorem:main} then
follows from the comparability of Theorem~\ref{theorem:combinatorics}.
We conclude with 
applications to the study of geometric limits,
proving Theorem~\ref{theorem:geometric:limit}.

\bold{Acknowledgements.}  I would like to thank Yair Minsky, in
particular, for introducing me to Thurston's conjecture and for
suggesting the use of curve hierarchies to improve the volume lower bound.  
Thanks also to Bill Thurston, Steve
Kerckhoff, and Curt McMullen for useful conversations, and to Lewis Bowen and
the referee for corrections and suggestions.

\section{The extended Weil-Petersson metric}
\label{section:preliminaries}
Let $S$ be a compact oriented surface of negative Euler
characteristic.  We allow $S$ to have boundary and let $\interior(S)$
denote its interior.  Let $\eS$ denote the set of isotopy classes of
essential, non-peripheral, simple closed curves on $S$.  

A {\em pants decomposition} $P
\subset \eS$ is a maximal collection of isotopy classes with pairwise
disjoint representatives on $S$.  The usual geometric intersection
number $i(\alpha,\beta)$ of a pair $(\alpha,\beta) \in \eS \times \eS$
generalizes to a total intersection number $i(P,P')$ of pants
decompositions by summing the geometric intersections of their
components.

The {\em Teichm\"uller space} $\Teich(S)$ of $S$ parameterizes finite
area hyperbolic structures on $S$ up to isotopy.
Points in $\Teich(S)$ are
pairs $(f,X)$  where $X$ is a finite area hyperbolic
surface $X$ equipped with a homeomorphism $f \colon \interior(S) \to X$, up
to the equivalence $(f,X) \sim (g,Y) $ if there is an isometry $\phi
\colon X \to Y$ for which $\phi \compos f \simeq g$.  
A pants decomposition $P = \alpha_1 \cup \ldots \cup \alpha_{|P|}$
determines {\em Fenchel-Nielsen coordinates} $$(\ell_X(\alpha_1),
\ldots, \ell_X(\alpha_{|P|}),\theta_X(\alpha_1),
\ldots, \theta_X(\alpha_{|P|})) \in \reals_+^{|P|} \times \reals^{|P|}$$
for each $X \in \Teich(S)$, indicating $X$ is assembled from
hyperbolic pairs of pants with boundary lengths prescribed by
$\ell_X(\alpha_i)$ glued together twisted by $\theta_X(\alpha_i)$.
(For more on Teichm\"uller space and Fenchel-Nielsen coordinates see
\cite{Imayoshi:Taniguchi:book} or \cite{Gardiner:book}).

\bold{The Weil-Petersson metric.}
Each $X \in \Teich(S)$ is naturally a complex 1-manifold via its
uniformization $X = \half^2/\Gamma$ as the quotient of the upper half
plane by a Fuchsian group.  The Teichm\"uller space has  a
complex manifold structure of dimension $3g-3+n$ where $S$ has genus
$g$ and $n$ boundary components.  

The space of holomorphic quadratic differentials 
$Q(X)$ on $X \in \Teich(S)$ (holomorphic forms of type $\phi(z) dz^2$
on $X$) is 
naturally the cotangent space $T^*_X \Teich(S)$ to $\Teich(S)$ at
$X$.    The {\em Weil-Petersson metric} on
$\Teich(S)$ unifies 
the hyperbolic and holomorphic perspectives on $X$:  
it arises from the $L^2$ inner product
on $Q(X)$, namely
$$\langle\varphi,\psi \rangle_{\rm WP} =  \int_X \frac{\varphi \bar{\psi}}{\rho^2}
$$ 
where $\rho(z) |dz|$ is the hyperbolic metric on $X$, by 
the usual pairing
$$(\mu,\varphi)_X =  \int_X \mu \varphi $$
between $T_X \Teich(S)$ and $T^*_X \Teich(S)$ (see,
e.g. \cite[Sec. 1]{Wolpert:Nielsen}).  In what follows, we will be  
interested only in the Riemannian part $g_{\rm WP}$ of the
Weil-Petersson metric, and its associated distance function $d_{\rm
WP}(.,.)$ on $\Teich(S)$.

The Weil-Petersson metric has negative sectional curvature
\cite{Tromba:sectional} \cite{Wolpert:sectional}, 
and the {\em modular group} $\Mod(S)$ (the group of isotopy classes of
orientation preserving homeomorphisms of $S$) acts by
isometries
of $g_{\rm WP}$.
Thus, $g_{\rm WP}$ descends to a metric on the Moduli space $\calM(S)
=
\Teich(S)/\Mod(S)$.

Work of S. Wolpert shows two important properties of the
Weil-Petersson metric we will use:

\bold{WPI} {\em The Weil-Petersson metric is not complete:}  
``pinching geodesics'' in the Teichm\"uller metric (which leave every
compact set of $\Teich(S)$) have finite Weil-Petersson length
\cite{Wolpert:noncompleteness}.

\bold{WPII} {\em The Weil-Petersson metric is geodesically convex:} in 
fact, for $\alpha \in \eS$ the length function
$\ell_{(.)}(\alpha)$ is strictly convex along Weil-Petersson geodesics
\cite{Wolpert:Nielsen}. 

\medskip

\bold{The augmented Teichm\"uller space.}  In \cite{Masur:WP},
H. Masur shows the Weil-Petersson metric
extends to the {\em augmented  
Teichm\"uller space} $\closure{\Teich(S)}$ obtained by adding boundary
Teichm\"uller spaces consisting of marked {\em noded} Riemann
surfaces, which we now describe
 (see \cite{Bers:nodes} for a detailed discussion).

A {\em Riemann surface with nodes} $W$ is a connected complex space so
that each point $p \in W$ has a neighborhood isomorphic to $\{ z \in
\cx \st |z| <1 \}$ or isomorphic to $\{ (z,w) \in \cx^2 \st |z|<1, |w|<1,\
\text{and} \ zw =0\}$ by an isomorphism sending $p$ to $(0,0) \in
\cx^2$.  In the latter case, $p$ is called a {\em node}  
of $X$.  The complement of the nodes is a union of Riemann surfaces
called the {\em pieces} of $W$.  We say $W$ is {\em hyperbolic} if
each piece of $W$ admits a complete finite-area hyperbolic structure.
 
Given curves $\nu_1,\ldots,\nu_j$ in a pants decomposition $P =
\{\alpha_1, \ldots, \alpha_{|P|}\}$ of $S$, 
a {\em marked noded hyperbolic surface pinched along
$\nu_1,\ldots,\nu_j$}
is a noded hyperbolic Riemann surface $W$ together with a continuous map
$$f \colon \interior(S) \to W$$
so that $f\vert_{S - \nu_1 \cup \ldots \cup \nu_j}$ is a
homeomorphism on to the union of the pieces of $W$.  
Let $S - \calN(\nu_1) \cup \ldots \cup \calN(\nu_j) = S_1
\cup \ldots \cup S_k$ where $\calN(\nu_i)$ are pairwise disjoint
open collars about each $\nu_i$.  Then the pair $(f,W)$ determines
a point $$\Teich(S_1) \times \ldots \times \Teich(S_k)$$ in the
product Teichm\"uller space by taking the restriction of $f$ to each
component of $\interior(S) - \nu_1 \cup \ldots \cup \nu_j$ as a
marking on each piece of $W$.  A marked piece of $W_l \in
\Teich(S_l)$, $1\le l \le k$, has Fenchel-Nielsen coordinates with
respect to the elements of the pants decomposition $P$ that lie in $S_l$.

Two marked hyperbolic noded surfaces
$(f_1,W_1)$ and $(f_2,W_2)$ are equivalent, if there is continuous map 
$\phi \colon W_1 \to W_2$ that is isometric on each piece of $W_1$ for 
which $\phi \compos f_1 = f_2$ after precomposition with an isotopy of
$S$.   

The {\em augmented Teichm\"uller space} $\closure{\Teich(S)}$ is
obtained by adjoining equivalence classes of marked noded hyperbolic
surfaces to $\Teich(S)$.  The topology on $\closure{\Teich(S)}$ is
given as follows.  Given a pants decomposition $P$, and a point $W \in
\closure{\Teich(S)}$ with curves $\nu_1 \cup \ldots \cup \nu_j$ 
in $P$ pinched to nodes, we extend the Fenchel-Nielsen coordinates to
$W$ by defining the coordinates $\ell_W(\nu_i) =
0$.  Then a neighborhood 
of $W$ in $\closure{\Teich(S)}$ consists of (possibly noded)
hyperbolic Riemann surfaces $X$ whose
length coordinates $\ell_X(\alpha_p)$ are close to those of $W$
for $p = 1,\ldots,3g-3+n$,
and whose twist coordinates $\theta_X(\alpha_p)$ are close to those of
$X$ for each $p$ such that  
$\alpha_p \not= \nu_i.$
(see \cite[App. B]{Imayoshi:Taniguchi:book}). 

\smallskip

The Weil-Petersson metric extends to the augmented Teichm\"uller space
as its completion (see \cite{Masur:WP}), giving a $\Mod(S)$ invariant
metric on $\closure{\Teich(S)}$.  The quotient
$$\closure{\Teich(S)}/\Mod(S) = \closure{\calM(S)},$$ the familiar
Mayer-Mumford-Deligne  
compactification of the moduli space (see \cite{Bers:nodes}), inherits
a complete extension of the Weil-Petersson metric on $\calM(S)$.  We
denote the corresponding distance by $$d_{\closure{\rm WP}} \colon
\closure{\Teich(S)} \times
\closure{\Teich(S)} \to \reals_{\ge 0}.$$

Evidently, the failure of completeness of the Weil-Petersson metric
occurs at limits of {\em pinching sequences} $X_t$ for which 
the length coordinates $\ell_{X_t}(\nu_i)$ tend to zero for some
collection of curves in a pants decomposition $P$. 

Given a pants decomposition $P$ and a collection $\alpha_1,
\ldots, \alpha_k$ of curves in $P$, the minimal distance from a point
$X \in \Teich(S)$ to a noded Riemann
surface $Z$ with nodes along $\alpha_1,\ldots,\alpha_k$ is estimated
in terms of the geodesic length sum 
\begin{equation*}
\ell = \ell_X(\alpha_1) + \ldots + \ell_X(\alpha_k)
\end{equation*}
of the lengths of $\alpha_i$ on $X$ by
\begin{equation}
d_{\overline{\rm WP}}(X,Z) = \sqrt{2\pi \ell} + O(\ell^2)
\label{equation:distance}
\end{equation}
(see \cite[Cor. 21]{Wolpert:compl}).

\bold{Remark:} This estimate is a recent improvement of similar estimates
originally obtained in \cite[Ex. 4.3]{Wolpert:universal} and cited in
earlier versions of this manuscript.

\bold{Sub level sets.}  We recall the following theorem of Bers.

\begin{theorem}[Bers]
There is a constant $L >0$ depending only on $S$ such that for any $X
\in \Teich(S)$ there is a pants decomposition $P$ such that 
$\ell_X(\alpha) < L$ for each $\alpha \in P$.
\label{theorem:bers:constant}
\end{theorem}
We call this $L$ the {\em Bers constant} for $S$.

Given a pants decomposition $P$, and a positive real number $\ell \in
\reals_+$, we consider the sub level set
$$V_\ell(P) = \left\{ X \in \Teich(S) \st \max_{\alpha \in P} \left\{
\ell_X(\alpha) \right\} < \ell \right\}.$$
Then by Bers's theorem, the union of the sets $V_L(P)$ over all pants
decompositions gives an open cover of $\Teich(S)$.  Because $L$
depends only on $S$, we abbreviate
$$V(P) = V_L(P).$$

Then we have the following:

\begin{prop} The sub level sets $V(P)$ have the following properties.
\begin{enumerate} 
\item Each $V(P)$ is convex in the Weil-Petersson metric, and
\item there is a constant $D>0$, depending only on $S$, for which 
the Weil-Petersson diameter $\diam_{\rm WP}(V(P)) < D.$ 
\end{enumerate}
\label{prop:convex:bounded}
\end{prop}
\bold{Proof:}
Geodesic convexity of $V(P)$ follows immediately from {\bf WPII}, the
convexity of the geodesic length functions $\ell_X(.)(\alpha)$ for each 
$\alpha \in P$.

To see each $V(P)$ has bounded Weil-Petersson diameter, let $W_P$ be
the (unique) maximally noded Riemann surface where each curve in $P$
is pinched.  By equation~\eqref{equation:distance} there is a
constant $C(L)$ so that
for each $X \in V(P)$ we have 
$$d_{\closure{\rm WP}}(X,W_P) < C(L).$$ 

By the triangle inequality for $d_{\closure{\rm WP}}$, if $X$ and $Y$
lie in $V(P)$ then the distance $d_{\closure{\rm WP}}(X,Y)$ is bounded
by $2C(L)$.  By geodesic convexity of $V(P)$ the geodesic joining
$X$ to $Y$ lies in $V(P)$, so we have the bound $$d_{\rm WP}(X,Y) < 2C
(L)$$ on $d_{\rm WP}(X,Y)$ which we set equal to $D$.

\qed

\section{A combinatorial Weil-Petersson distance}
\label{section:combinatorics} 
In this section, we relate the coarse geometry of the Weil-Petersson
metric to the {\em pants graph} ${\bf P}(S)$ defined in the
introduction.  We do this by exhibiting a {\em quasi-isometry} between 
the two spaces with their respective distances.

\begin{defn}
Given $k_1 >1$ and $k_2 >0$, a map $f \colon (X,d) \to (Y,d')$ of
metric spaces is a $(k_1, k_2)$-{\em 
quasi-isometric embedding} if for each pair of points $x$ and $y$ in
$X$ we have 
$$\frac{d(x,y)}{k_1} - k_2 \le d'(f(x),f(y)) \le k_1 d(x,y) + k_2.$$
\end{defn}

The spaces $(X,d)$ and $(Y,d')$ are {\em quasi-isometric} if for some
$k_1 >1$ and $k_2>0$ there are $(k_1,k_2)$-quasi-isometric embeddings from
$(X,d)$ to $(Y,d')$ and from $(Y,d')$ to $(X,d)$.  In practice, it
suffices to exhibit a {\em quasi-isometry} from $(X,d)$ to
$(Y,d')$, namely, a quasi-isometric embedding with uniformly dense
image.  Given such a quasi-isometry from $(X,d)$ to  
$(Y,d')$, a quasi-isometric embedding from $(Y,d')$ to $(X,d)$ is
readily constructed, so the spaces are quasi-isometric.
\medskip

Let $$Q \colon {\bf P}^0(S) \to \Teich(S)$$ be any embedding of
the vertices ${\bf P}^0(S)$ of ${\bf P}(S)$ into $\Teich(S)$
so that $Q(P)$ lies in $V(P)$.  The main theorem of this section is
the following:
\begin{theorem}
The map $Q$ is a quasi-isometry of the $0$-skeleton  ${\bf P}^0(S)$
of ${\bf P}(S)$ with $\Teich(S)$ with its Weil-Petersson distance.
\label{theorem:pants:qi}
\end{theorem}

\bold{Proof:}
By the uniform bound 
$$\diam_{\rm WP}(V(P))< D$$ 
on the diameter of $V(P)$, the image $Q({\bf P}(S))$ 
is $D$-dense in $\Teich(S)$.  It suffices, then, to show
that there are uniform constants $A_1 \ge 1$ and $A_2 \ge 0$ so that
$$\frac{1}{A_1}d_{\bf P}(P_0,P_1) - A_2 \le d_{\rm WP}(Q(P_0),Q(P_1)) \le
A_1 d_{\bf P}(P_0,P_1) + A_2.$$

We first show that the map $Q$ is $2D$-Lipschitz.  Given $P_0$ and
$P_1$ such that 
$$d_{\bf P}(P_0,P_1) = 1,$$ 
$P_0$ and $P_1$ differ 
by a single elementary move.  Let $\alpha \in P_0$ and $\beta \in
P_1$ be the curves involved in this elementary move, i.e.
$P_0 - \alpha = P_1 - \beta$ and $i(\alpha,\beta) = 1$ or $2$
depending on whether the component $S_\alpha \subset S - (P_0 -
\alpha)$ containing $\alpha$ is a punctured torus or
four-times punctured sphere.  

Let $Z \in \Teich(S_\alpha)$ be the ``square'' punctured torus: i.e.
$Z$ is obtained by identifying opposite sides of an ideal square in
$\half^2$ with order-$4$ rotational symmetry about the origin in the
disk model of $\half^2$.  
\makefig{The cover $\wt Z$ with lifts of $\alpha$ and $\beta$ passing through
the origin.}{figure:square}{\psfig{file=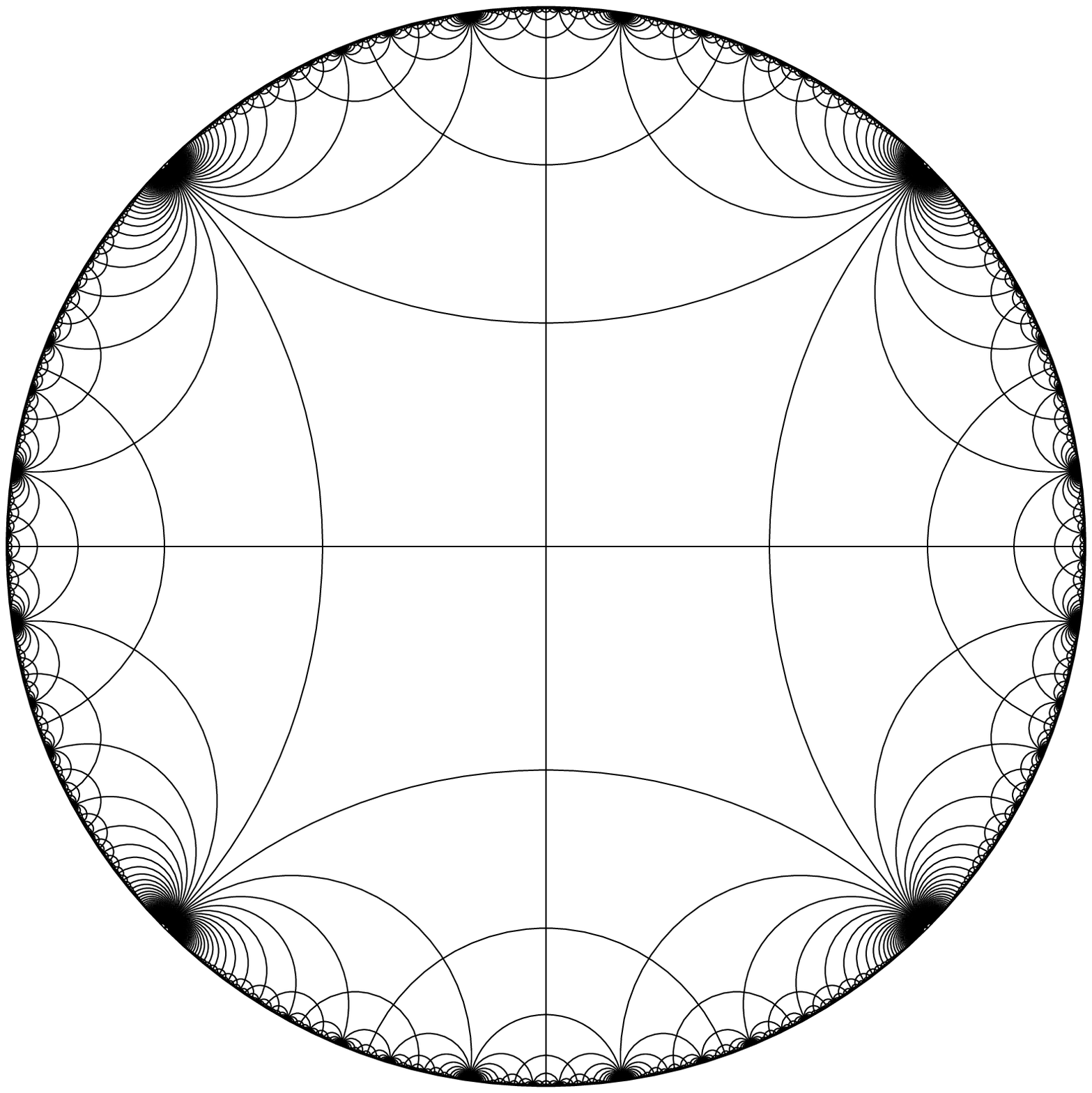,height=3in}}
Marking $Z$ so that the common
perpendiculars to the opposite sides descend to closed geodesics
$\alpha$ and $\beta$ on $Z$ (see Figure~\ref{figure:square}), we have by symmetry that the length
$\ell_Z(\alpha)$ equals the length $\ell_Z(\beta)$, which is the
shortest length of any non-peripheral simple closed curve on $Z$.  In
particular, we have 
$$\ell_Z(\alpha)< L(S_\alpha)$$ 
where $L(S_\alpha)$ is the Bers constant for $S_\alpha$.

Let $(f, W) \in
\overline{\Teich(S)}$ be the noded Riemann surface with nodes at each
$\gamma \in P - \alpha$, and one non-rigid piece $f \vert_{S_\alpha}
\to Z$.
For any $Z' \in \Teich(S)$ sufficiently close to $Z$, if $Z' \in
V(P')$ then $P'$ contains an essential non-peripheral curve in
$S_\alpha$.  Taking $Z'$ arbitrarily close to $Z$, we may conclude
that $L(S_\alpha) \le L(S) = L$.  Then for any Riemann surface $Z' \in
\Teich(S)$ sufficiently close to $W \in \overline{\Teich(S)}$ we have
$$\ell_{Z'}(\gamma) < L$$ for each $\gamma \in P_0 \cup P_1$.  In
other words, $Z'$ lies in the intersection $V(P_0) \cap V(P_1)$.
Letting $Z$ be the double of the symmetric ideal square described
above, a similar argument handles the genus-0 case.

Therefore we may conclude that 
$$d_{\rm WP}(Q(P_0),Q(P_1)) <2D \ \ \ \text{ when } \ \ \ d_{\bf P}(P_0,P_1)
=1,$$ so by the triangle inequality, $Q$ is $2D$-Lipschitz.
To show that for some $A_1$ and $A_2$ the inequality
$$\frac{1}{A_1}d_{\bf P}(P_0,P_1) - A_2\le d_{\rm WP}(Q(P_0),Q(P_1))$$
holds is somewhat more delicate.  We break this into a series of
lemmas.  

\begin{lem}
Given $L' > L$,
there is an integer $B>0$ so that given $P$ and $P'$ in
${\bf P}(S)$ for which $V_{L'}(P)\cap V_{L'}(P') \not= \nullset$, we have
$d_{\bf P}(P,P') \le B.$
\label{lemma:unit:length}
\end{lem}

\bold{Proof:}  
The hypotheses imply that there is some $X \in \Teich(S)$ so that
$\ell_X(\alpha) < L'$ for each $\alpha \in P \cup P'$.
By an application of the collar lemma
\cite[Thm. 4.4.6]{Buser:book:spectra} there is a constant $C$
depending only $L'$ and $S$ so that the total geometric 
intersection number $i(P,P')$ satisfies
$$i(P,P') \le C.$$

Let ${\bf Tw}(P) \cong \zed^{|P|}$ denote the subgroup of $\Mod(S)$
generated by Dehn twists about the curves in $P$.  
Then the function $i(P, .) \colon {\bf P}^0(S) \to \zed$
descends to a function
$$i(P, .) \colon {\bf P}^0(S)/{\bf Tw}(P) \to \zed$$
whose sub level sets are bounded: in other words
there are only finitely many equivalence classes $$\{[P_1],\ldots,
[P_c]\} \subset {\bf P}^0(S)/{\bf Tw}(P)$$ for which
$i(P,[P_j]) \le C.$

Since $d_{\bf P}(P,.)$ also descends to a function
$$d_{\bf P}(P,.) \colon {\bf P}^0(S)/{\bf Tw}(P) \to \zed$$
on ${\bf P}^0(S)/{\bf Tw}(P)$, we have
$$d_{\bf P}(P,P') \le B$$ where
$$B = \max_{j=1,\ldots,c}\{d_{\bf P}(P,[P_j])\}.$$
\qed

\begin{lem}
Given $L' > L$, there is an integer $J>0$, so that if
$X_t$, $t \in [0,1]$, is a unit-length Weil-Petersson geodesic
joining $X_0$ and $X_1$, then there exist pants decompositions $P_1,
\ldots, P_J$ so that $\{X_t\}_{t=0}^1$ lies in the union 
$$V_{L'}(P_1) \cup
\ldots \cup V_{L'}(P_J).$$
\label{lemma:bounded:covering}
\end{lem}

\bold{Proof:}  
Recall from Theorem~\ref{theorem:bers:constant} that the sets $V(P)
\subset V_{L'}(P)$  
cover $\Teich(S)$.
Let $P_1, \ldots, P_m$ determine sets $V(P_1), \ldots, V(P_m)$
so that for each $t \in [0,1]$ we have
$$X_t \in V(P_1) \cup \ldots \cup V(P_m).$$

Let $$d_{L',P} \colon \bdry V(P) \to \reals_+$$
be the function 
$$d_{L',P}(X) = \inf_{Y \in \bdry V_{L'}(P)}
d_{\rm WP}(X, Y).$$

We claim that there is an $\epsilon_0$ depending only on $L$, and $L'$
so that $$d_{L',P}(X) >\epsilon_0.$$ The function $d_{L',P}(X)$
naturally extends to the metric completion $\closure{\bdry V(P)}$ of
$\bdry V(P)$, and $d_{L',P}(X)$ is invariant under the action of ${\bf
Tw}(P)$.

Let $\{(\ell_i,\theta_i) \in \reals^{|P|}_+ \times \reals^{|P|}\}$
denote Fenchel-Nielsen coordinates for $\Teich(S)$ adapted to the
pants decomposition $P$.  To extend these Fenchel Nielsen coordinates
to the completion, we denote by 
$$\reals_{\ge 0} \times \reals /\sim$$ the quotient
of $\reals_{\ge 0} \times \reals$ by the equivalence relation
$(0,\theta) \sim (0, \theta')$.
Then the completion $\closure{V_{L'}(P)}$ of
$V_{L'}(P)$ in $\closure{\Teich(S)}$ admits {\em extended}
Fenchel-Nielsen coordinates 
$$\closure{V_{L'}(P)} = \{ (\ell_i, 
\theta_i) \in \reals^+ \times \reals/\sim \st \ell_i \le L', \ i =
1,\ldots, |P|\}$$ 
where each point with $\ell_j = 0$ for some $j$ lies
in the completion.  The extended isometric action of ${\bf Tw}(P)$ on
$\closure{\Teich(S)}$ is cocompact on $
\closure{V_{L'}(P)}$, since ${\bf Tw}(P)$ preserves each length
coordinate and acts by translations on each twist coordinate.

In these extended Fenchel-Nielsen coordinates, the completion
$\closure{\bdry V_{L'}(P)}$ of the boundary $\bdry V_{L'}(P)$ is the
locus of coordinates for which $\ell_j = L'$ for some $j \in
1,\ldots,|P|$.  Thus $$\closure{\bdry V_{L'}(S)}/{\bf Tw}(P)$$ is a
closed subset of the compact set $\closure{V_{L'}(P)}/{\bf Tw}(P)$ and
is thus compact.

Since the quotients
$$\closure{\bdry V_{L'}(P)}/{\bf Tw}(P) \ \ \ \text{and} \ \ \ 
\closure{\bdry V(P)}/{\bf Tw}(P)$$
are disjoint compact subsets of $\closure{V_{L'}(P)}/{\bf Tw}(P)$,
it follows that  the 
function 
$$\calI_{L'} \colon {\bf P}^0(S) \to \reals_+$$ given by
$$\calI_{L'}(P) = \inf_{X \in {\bdry V(P)}} d_{L',P}(X)$$
is positive.
But $\calI_{L'}$ is $\Mod(S)$-invariant, 
so it descends to a function 
$$\calI_{L'} \colon {\bf P}^0(S) / \Mod(S) \to \reals_+.$$
Since 
${\bf P}^0(S)/\Mod(S)$ is finite,
we may set $\epsilon_0$ equal to the
infimum of $\calI_{L'}([P])$ over the finite number of equivalence
classes $[P] \in {\bf P}^0(S)/\Mod(S)$.

If $X_{t_0}$ lies in $V(P)$, then, $X_t$ lies in $V_{L'}(P)$, provided
$t$ lies in $(t_0 -\epsilon_0, t_0 + \epsilon_0)$.  It follows that
after setting $J$ equal to the least integer greater than $2/\epsilon_0$, we
may select from the pants decompositions $P_1, \ldots, P_m$ pants
decompositions $P_1, \ldots, P_J$ (possibly with repetition) so that
$$X_t \in V_{L'}(P_1) \cup \ldots \cup V_{L'}(P_J)$$ for each $t \in [0,1]$.

\qed

To complete the proof of Theorem~\ref{theorem:pants:qi}, let $X_t$ be
the Weil-Petersson geodesic joining 
arbitrary distinct Riemann surfaces $X$ and $Y$ in $\Teich(S)$.  Let
$P_X$ and $P_Y$ be pants decompositions for which $X$ lies in
$V(P_X)$ and $Y$ lies in $V(P_Y)$.  Let $I(P) \subset [0,1]$
denote the values of $t$ for which $X_t
\in V_{2L}(P)$.  By convexity of $V(P)$
(Proposition~\ref{prop:convex:bounded}) each  
$I(P)$ is an interval.  

Taking $L' = 2L$, Lemma~\ref{lemma:bounded:covering} provides a $J>0$
and a sequence $\{P_j \}_{j=0}^N \in {\bf P}(S)$ so that 
\begin{itemize}
\item $X(t)$ is covered by the union $\cup_j I(P_j)$,
\item the least upper bound of $I(P_j)$ lies in $I(P_{j+1})$, and
\item $N \le J(d_{\rm WP}(X,Y) +1)$.
\end{itemize}
Thus, we have
\begin{equation}
\frac{N}{J} -1 \le d_{\rm WP}(X,Y).
\label{equation:covering}
\end{equation}
Moreover, for successive pants decompositions $P_j, P_{j+1}$,
we have $$V_{2L}(P_j) \cap V_{2L}(P_{j+1}) \not= \nullset,$$
so applying Lemma~\ref{lemma:unit:length} with $L' = 2L$, we have a $B 
>0$
for which
\begin{equation}
d_{\bf P}(P_X,P_Y) \le B N.
\label{equation:partition}
\end{equation}
Combining equations \ref{equation:covering}
and \ref{equation:partition}, we have
\begin{equation}
\frac{d_{\bf P}(P_X,P_Y)}{B J} - 1 \le d_{\rm WP}(Q(P_X),Q(P_Y)),
\end{equation}
where $B$ and $J$ depend only on $L$ which depends only on $S$.
Setting $A_1 = B J $ and $A_2 = 1$ concludes the proof of Theorem~\ref{theorem:pants:qi}.
\qed

\section{Bounding the core volume from below}
\label{section:lower}
To simplify notation, let $$C(X,Y) = \core(Q(X,Y))$$ and recall that
 $\vol(X,Y)$ denotes the convex core volume $
\vol(\core(Q(X,Y))).$ 
In this section we prove
\begin{theorem}
Given $S$, there are constants $K_1 >1$ and $K_2 >0$ so that
$$\frac{1}{K_1} d_{\rm WP}(X,Y) - K_2 \le \vol(X,Y).$$
\label{theorem:volume:lower}
\end{theorem}
The proof is given as a series of lemmas.

Fix attention on a given quasi-Fuchsian manifold $Q(X,Y)$.  Given a
constant $L_0 >0$, let $\eS_{< L_0} \subset \eS$ 
denote the set of isotopy classes
$$\eS_{< L_0} = \{ \alpha \in \eS \st \ell_{Q(X,Y)}(\alpha) < L_0 \}.$$
Theorem~\ref{theorem:volume:lower} will follow from a linear lower bound on
$\vol(X,Y)$  given in terms of the size of $\eS_{<L}$
(Lemma~\ref{lemma:isotopy:classes}) and the following lemma.

\begin{lem}
Let $P_X$ and $P_Y$ be pants decompositions so that $X \in V(P_X)$ and
$Y \in V(P_Y)$.   Then there is a constant $K$ depending
only on $S$ so that 
$$d_{\bf P}(P_X,P_Y) \le K |\eS_{< L}|.$$ 
\label{lemma:count}
\end{lem}

The lemma will follow from the following general result on paths in
${\bf P}(S)$ that are built out of a given collection of curves in
$\eS$.
For reference, let $$\pi_\eS
\colon {\bf P}^0(S) \to \eS$$ denote the projection that assigns to
each $P$ the collection of curves used to build it.

\begin{lem}
Let $k \in \natls$ and let $g = \{ P_I = P_0,\ldots , P_N = P_T \} \subset {\bf P}(S)$ be a
sequence of pants decompositions 
with the property that $d_{\bf P}(P_j,P_{j+1}) <  k$.  Let $\eS_g
\subset \eS$ denote the image of $g$ under the 
projection $\pi_\eS$.
There is a constant $K_0 >0$ depending on $k$ and $S$ so that
$$d_{\bf P}(P_I,P_T) < K_0 |\eS_g|.$$
\label{lemma:projection:argument}
\end{lem}

\bold{The complex of curves.}  
The graph ${\bf P}(S)$ is related to the {\em complex of curves}
$\calC(S)$, introduced by W. Harvey \cite{Harvey:CC}.  To prove
Lemma~\ref{lemma:projection:argument} we describe recent work of Masur
and Minksy on $\calC(S)$.  The main result of 
\cite{Masur:Minsky:CCI} shows that $\calC(S)$ is in fact a {\em Gromov
hyperbolic} metric space with the metric obtained by making each
simplex a standard Euclidean simplex.  Its sequel
\cite{Masur:Minsky:CCII} 
introduces a theory of
{\em hierarchies} 
of so-called ``tight geodesics'' in
$\calC(S)$ and in sub-complexes $\calC(Y)$ for essential subsurfaces $Y
\subset S$.  Such hierarchies and their hyperbolicity properties play
an integral role in our control of volume.

To describe the topological type of $S$, we let $$d(S) =
\dim_\cx(\Teich(S)) = 3g-3 + n$$ where $S$ has genus $g$ with $n$
boundary components.  We consider only those surfaces $S$ for which
$\interior(S)$ admits a hyperbolic structure (so $d(S) >0$).  

The {\em complex of curves} $\calC(S)$ is a simplicial
complex with $0$-skeleton $\eS$, and higher dimensional simplices
described as follows:
\begin{itemize}
\item for $d(S) > 1$ and $k \ge 1$, $k$-simplices of $\calC(S)$ span
$k+1$-tuples 
$\alpha_1, \ldots, \alpha_{k+1}$ of vertices for which
$i(\alpha_i,\alpha_j) = 0$, and
\item if $d(S) = 1$, $\calC(S)$ is a 1-complex whose edges join vertices
$\alpha$ and $\alpha'$ in $\calC(S)$ that intersect minimally;
i.e. $\calC(S) = {\bf P}(S)$.
\end{itemize}

Given an essential subsurface $Y \subset S$ with $d(Y) \ge 2$, the curve
complex $\calC(Y)$ is naturally a subcomplex of $\calC(S)$.  Given a
set $W$ let $\calP(W)$ denote its power set, i.e. the set of all
subsets of $W$.  Masur and Minsky
define a projection $$\pi_Y \colon \calC(S) \to \calP(\calC(Y)),$$ by
setting $\pi(\alpha) = \alpha$ if $\alpha \in \calC(Y)$, and taking
$$\pi_Y(\alpha) = \bigcup_{\alpha' \subset \alpha \cap Y} \bdry
\calN(\alpha' \cup \bdry_{\alpha'} Y)$$
where $\alpha'$ is an arc of essential intersection of $\alpha$ with
$Y$, $\bdry_{\alpha'} Y \subset \bdry Y$ is the components of $\bdry
Y$ that $\alpha'$ meets in its endpoints, and $\calN(.)$ denotes a
regular neighborhood of their union (see
\cite[Sec. 2]{Masur:Minsky:CCII}).

When $A$ and $B$ are two subsets of $\calC(Y)$,
\cite{Masur:Minsky:CCII} defines a coarse
distance $d_Y(A,B)$ by taking the diameter $$d_Y(A,B) =
\diam_{\calC(Y)}(A \cup B),$$ 
in $\calC(Y)$ of $A$ and $B$.
Note that while $d_Y(.,.)$
is more a diameter than a distance when $A$ and $B$ are close, 
it gives a useful notion of distance between sets of bounded diameter
and does satisfy the triangle inequality.

By \cite[Lem. 2.3]{Masur:Minsky:CCII}
the projection $\pi_Y$ has a Lipschitz property: 
if $\Delta$ is a
simplex in $\calC(S)$ so that $\Delta$ intersects $Y$, then we have
$\diam_{\calC(Y)}(\pi_Y(\Delta)) \le 2$.  
If $P_I$ and $P_T$ are two subsets of $\calC(S)$, letting $\pi_Y(P_I) = 
\cup_{\alpha \in P_I} \pi_Y(\alpha)$, and likewise for $P_T$,
then the {\em projection distance} $d_Y(P_I,P_T)$ between $P_I$ and
$P_T$ (or {\em distance in $Y$}) 
is defined by $$d_Y(P_I,P_T) =
d_Y(\pi_Y(P_I) , \pi_Y(P_T)).$$ In particular, if $P$ and $P'$ are
pants decompositions that differ by a single elementary move, then we
have $d_Y(P,P') \le 4$ (see \cite[Lem. 2.5]{Masur:Minsky:CCII}).

A central theorem we will use is the following:
\begin{theorem}[Thm. 6.12 of \cite{Masur:Minsky:CCII}]
There is a constant $M_0(S)$ so that given $M> M_0$ there exist $c_0$
and $c_1$ so that if $P_I$ and $P_T$ are pants decompositions in
${\bf P}(S)$ then we have $$\frac{1}{c_0}d_{\bf P}(P_I,P_T) - c_1
\le
\sum_{\stackrel{Y \subseteq S}{d_Y(P_I,P_T) >M}} d_Y(P_I,P_T) \le
c_0 d_{\bf P}(P_I,P_T) + c_1$$ where the sum is taken over all
non-annular essential subsurfaces $Y \subseteq S$ satisfying
$d_Y(P_I, P_T) > M$.
\label{theorem:MMII:bound}
\end{theorem}

We apply this result to prove Lemma~\ref{lemma:projection:argument}.
Our argument is quite similar to that of
\cite[Thm. 6.10]{Masur:Minsky:CCII}, where it is shown that a given
pants decomposition along an elementary move sequence can contribute
to progress in only boundedly many projections to subsurfaces
simultaneously.  We seek the analogous statement for a single curve
occuring in pants decompositions joining $P_I$ to $P_T$.

\bold{Proof:} {\em (of Lemma~\ref{lemma:projection:argument})}.
To prove the lemma, we will relate the sum of the projections to the
size of $\eS_g$.  To do this, we note that when the projection distance
$d_Y(P_I,P_T)$ is large, there must be a definite portion of the
projection of $g$ to $Y$ that is far from both $\pi_Y(P_I)$ and
$\pi_Y(P_T)$ in $d_Y(.,.)$; this follows from the triangle inequality
for $d_Y$ and the fact that elementary moves in ${\bf P}(S)$ make
Lipschitz progress as measured by $d_Y$.

We argue that a given curve $\alpha$ can contribute only to a bounded
amount of progress in boundedly many different subsurfaces.  Precisely, 
let $Y$ and $Z$ be two
essential, intersecting, non-annular subsurfaces of $S$, neither of
which is contained in the other.  
A lemma of Masur and Minsky \cite[Lem. 6.11]{Masur:Minsky:CCII}
enforces a partial ordering ``$\prec$'' on such subsurfaces with respect
to the pants decompositions $P$ and $P'$, provided the projection
distances
$d_Y(P,P')$ and $d_Z(P,P')$ are greater than a constant $M_2$
depending only on $S$.  
Taking $M_2$ to be the constant of
\cite[Lem. 6.2]{Masur:Minsky:CCII} with the same name, we say the
subsurfaces $Y$ and $Z$ are
{\em $(P,P')$-ordered} if we have $$d_Y(P,P') > M_2\ \ \ \text{and} \
\ \ d_Z(P,P') >M_2.$$

We rephrase \cite[Lem. 6.11]{Masur:Minsky:CCII} as follows.

\begin{lem}
There is a constant $M_3$ depending only on $S$ so that if $Y$ and $Z$
are $(P,P')$-ordered then one of two cases obtains.  Either $Y \prec
Z$, and we have $$d_Y(\bdry Z , P') < M_3 \ \ \text{and} \
\ d_Z(P, \bdry Y) < M_3,$$ or
$Z \prec Y$, and we have
$$d_Z(\bdry Y , P') < M_3 \ \ \text{and} \ \ d_Y(P, \bdry Z) <
M_3.$$

\label{lemma:order}
\end{lem}

\bold{Proof:} 
If $Y$ and $Z$ are $(P,P')$-ordered in the above sense,
then by \cite[Lem. 6.2]{Masur:Minsky:CCII} they appear as {\em
domains} $Y = D(h)$ and $Z = D(k)$ supporting {\em tight geodesics} $h
\subset \calC(Y)$ and $k \subset \calC(Z)$ in any {\em hierarchy} $H$
(without annuli \cite[Sec. 8]{Masur:Minsky:CCII}) joining $P = I(H)$
to $P' = T(H)$.  The condition 
that $Y$ and $Z$ intersect and are non-nested guarantees that $h$ and
$k$ are {\em time-ordered}
\cite[Lem. 4.18]{Masur:Minsky:CCII}  (in the sense of
\cite[Defn. 4.16]{Masur:Minsky:CCII}).  The lemma then follows from an
application of \cite[Lem. 6.11]{Masur:Minsky:CCII} where $Y \prec Z$
represents the case $h \prec_t k$ and $Z \prec Y$ represents the case
$k \prec_t h$. 
\qed

Let $M_4 = M_2 + 2M_3 +4$, and let $M =\max\{4M_4, M_0\}$.
Consider an essential subsurface $Y \subseteq S$ for which
$d_Y(P_I,P_T) > M$.  As in the proof of
\cite[Lem. 6.10]{Masur:Minsky:CCII} let $J_Y$ denote the subset of
$[1,N]$ for which 
if $i \in J_Y$ then $P_i$ is ``deep'' in the projection to
$Y$: i.e. 
$$d_Y(P_I,P_i) > M_4 \
\ \text{and} \ \ d_Y(P_i,P_T) >M_4.$$ 


Given a subset $A \subset [1,N]$, we denote by 
$$\|A \|_Y = \diam_{\calC(Y)}(\{ \pi_Y(P_i) \st i \in A \})$$ 
the diameter of the projection of the pants decompositions with
indices in $A$ to the curve complex $\calC(Y)$.

Given $\alpha \in \eS_g$ for which $\pi_Y(\alpha) \not= \nullset$, we
denote by $J_Y(\alpha) \subset J_Y$ the subset for which if $i \in
J_Y(\alpha)$ then $\alpha$ lies in $P_i.$ 
We make three observations for later reference:
\begin{enumerate}
\item[{\bf I.}] By the Lipschitz property for $\pi_Y$, we have
$\|J_Y(\alpha)\|_Y \le 4.$ 
\item[{\bf II.}] If $i$ lies in $J_Y$ then there is some $\alpha \in \eS_g$
so that $i \in J_Y(\alpha)$.
\item[{\bf III.}] Since $d_{\bf P}(P_j,P_{j+1}) < k$, we have
$d_Y(P_j,P_{j+1}) < 4k$.
\end{enumerate}

Let $Y \subseteq S$ and $Z \subseteq S$ be two non-annular
intersecting subsurfaces neither of which is contained in the other
so that each contributes to the sum of Theorem~\ref{theorem:MMII:bound}:
i.e. we have
$$d_Y(P_I,P_T) >M\ \ \ \text{and} \ \ \ d_Z(P_I,P_T) >M.$$
This assumption guarantees, in particular, that 
$Y$ and $Z$ are $(P_I,P_T)$-ordered.
We make the following claim:
\begin{itemize}
\item[$(*)$]{\em If $J_Y(\alpha)$ is non-empty,
then $J_Z(\alpha)$ must be empty.}
\end{itemize} 
Arguing by contradiction, assume 
$$J_Y(\alpha) \not= \nullset \not= J_Z(\alpha).$$
We have 
$$\pi_Y(\alpha) \not= \nullset \not= \pi_Z(\alpha)$$  
so if $i \in J_Y(\alpha)$ and $j
\in J_Z(\alpha)$ then 
$\diam_Y(\pi_Y(P_i))\le2$, $\diam_Y(\pi_Y(P_j)) \le 2$ 
so it follows that
$d_Y(P_i,P_j) \le 4$,
since $\alpha$ lies in $P_i$ and in $P_j$.
The same conclusion holds with $Z$ in place of $Y$.

Since $i$ lies in $J_Y(\alpha)$ we have $d_Y(P_I,P_i) > M_4$ and
$d_Y(P_i,P_T) > M_4$, so it follows that
$$d_Y(P_I,P_j) \ge M_4 - 4 \ \ \ \text{and} \ \ \
d_Y(P_j,P_T) \ge M_4 - 4.$$ 
As $j$ lies in $J_Z(\alpha)$ we have
$$d_Z(P_I,P_j) \ge M_4,$$ 
and since $M_4 - 4 > M_2$, it follows that
$Y$ and $Z$ are also $(P_I,P_j)$-ordered.  Let $\prec_j$ denote the
$(P_I,P_j)$-ordering and assume without loss of generality that $Y
\prec_j Z$.  Then applying Lemma~\ref{lemma:order} we have $$d_Y(\bdry
Z, P_j) < M_3.$$

Since $Y$ and $Z$ are also
$(P_I,P_T)$-ordered, we may first assume that 
$Y\prec Z$.  Then Lemma~\ref{lemma:order} gives $d_Y(\bdry Z, P_T) <
M_3$ which implies that 
$$d_Y(P_j,P_T) < 2M_3,$$ contradicting the
assumption that $j \in J_Z(\alpha)$.  If on the other hand we have $Z
\prec Y$, then Lemma~\ref{lemma:order} gives $d_Y(P_I,\bdry Z) < M_3$
from which we conclude $$d_Y(P_I,P_j) < 2M_3,$$ which contradicts the
same assumption.
Thus, either $J_Y(\alpha)$ or $J_Z(\alpha)$ must be empty, and the
claim $(*)$ is proven.

\smallskip

Applying observations {\bf (I)} and {\bf (II)} above, 
if $\eS_Y =
\{\alpha \in \eS_g \st J_Y(\alpha) \not= \nullset \}$, then
we have 
$$\|J_Y\|_Y \le 4|\eS_Y|.$$ There is a uniform bound $s$ depending
only on $S$ to the size of any collection of subsurfaces any pair of
which is disjoint or nested (see \cite[Lem. 6.10,
proof]{Masur:Minsky:CCII}) so by our claim $(*)$, the number of $Y$
for which $J_Y(\alpha)$ can be non-empty is bounded by $s$.  Thus we
have $$\sum_{Y \subseteq S, J_Y \not= \nullset} \|J_Y\|_Y \le 4 s
|\eS_g|.$$

Applying observation {\bf (III)}, 
we have $d_Y(P_j,P_{j+1}) < 4k$.
Thus, for each $Y$ satisfying
$$d_Y(P_I,P_T) > M$$ the set $J_Y$ is in particular non-empty, and we have
$$d_Y(P_I,P_T) - 2 M_4 \le 4k\|J_Y\|_Y.$$
But since $4M_4 < M$ we have $2M_4 \le 4k\|J_Y\|_Y  $ and thus 
$$d_Y(P_I,P_T) \le 8k\|J_Y\|_Y.$$

Since $M = \max\{M_0,4M_4\}$, applying 
Theorem~\ref{theorem:MMII:bound} there are constants $c_0$ and $c_1$
so that we have
$$
\frac{1}{c_0}d_{\bf P}(P_I,P_T) - c_1 \le
\sum_{\stackrel{Y \subseteq S}{d_Y(P_I,P_T) > M}} d_Y(P_I,P_T) \le
32sk|\eS_g|.$$
Since $|\eS_g|$ is always at least $d(S)$, we may combine all of the
above constants into a single $K_0$ for which 
$$d_{\bf P}(P_I,P_T) \le K_0 |\eS_g|.$$

\qed

To prove Lemma~\ref{lemma:count}, we will apply
Lemma~\ref{lemma:projection:argument} to a sequence of pants
decompositions $\{P_X = P_0,\ldots,P_N = P_Y\}$ so that $\pi_\eS(P_j)
\subset \eS_{<L}$ for each $j$, and so that $d_{\bf P}(P_j,P_{j+1})$ is
bounded by an a priori constant.  The existence of such a sequence is
provided by an interpolation of $1$-Lipschitz homotopy equivalences of
hyperbolic surfaces into $Q(X,Y)$ that pass from one side of the
convex core to the other.  The existence of such an interpolation
follows from work of R. Canary on {\em simplicial hyperbolic surfaces}
which we now describe.

\bold{Simplicial hyperbolic surfaces.}
Let $\Sing_k(S)$ denote the finite-area marked {\em singular}
hyperbolic structures on $S$: complete finite area hyperbolic surfaces 
$Z$ with at most $k$ cone singularities, each with cone-angle at least
$2\pi$, equipped with marking homeomorphisms $h \colon \interior(S)
\to Z$ up to marking-preserving isometry.  
Roughly speaking, a {\em simplicial hyperbolic surface}
is a path-isometric mapping from a singular hyperbolic surface to a hyperbolic 3-manifold that is totally geodesic in the complement of a ``triangulation.'' 
We now make this notion precise.

Let $V = \{ v_1, \ldots, v_p \}$ be a finite subset of $S$.  Following
Hatcher \cite{Hatcher:triangulations}, an {\em essential arc} $\alpha$ 
in $(S,V)$ is an embedded arc meeting $\bdry S \cup V$ only in its
endpoints, which lie in $V$.   A collection $\{\alpha_0, \ldots,
\alpha_k\}$ of essential arcs in $(S,V)$ that are pairwise embedded and
non-isotopic {\em rel}-endpoints is called a {\em curve system}.
Let $\calA(S,V)$ denote the simplicial complex whose $k$-simplices 
$[\alpha_0, \ldots, \alpha_k]$ are 
curve systems with curves $\{\alpha_0,\ldots , \alpha_k\}$
with faces given by 
$k$-tuples of curves in $\{\alpha_0, \ldots, \alpha_k\}$.

If $V$ contains a point in each boundary component of the compact
surface $S$, then a {\em triangulation} of $S$ is
a maximal curve system in $\calA(S,V)$.  Likewise, we may view the interior
$\interior(S)$ of $S$ as a ``punctured surface''
by collapsing each boundary component $\gamma \subset
\bdry S$ to a point  $v_\gamma$ to obtain a surface $R$.
Then we have 
$$\interior(S) = R - \{v_\gamma \st \gamma \subset \bdry S\}.$$
If $V$ is a subset of $R$ containing $\cup_\gamma v_\gamma$, then a
{\em triangulation} of $\interior(S)$ is the restriction of a maximal
curve system in $\calA(R,V)$ restricted to $R - \{v_\gamma \st \gamma
\subset \bdry S\} = \interior(S)$.  
Note that in each definition, an edge may have its boundary vertices
identified and a face may have boundary edges identified.  

The main result of \cite{Hatcher:triangulations} guarantees that any
two triangulations in $\calA(S,V)$ are related by a finite sequence of
elementary moves (see Figure~\ref{figure:triangular:moves}).
\makefig{Elementary moves on triangulations.}{figure:triangular:moves}{\psfig{file=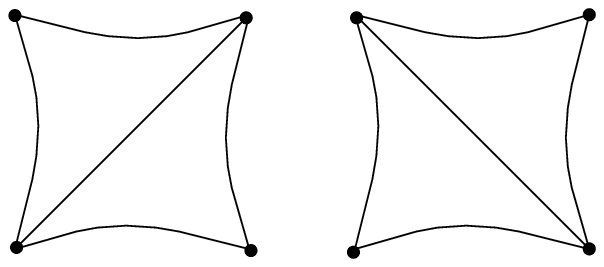,height=1.5in}} 

\smallskip
Let $T$ be a triangulation of $\interior(S)$ in the above sense, and
let $h \colon \interior(S) \to Z$ be a singular hyperbolic surface for
which $h$ is isotopic to a map with the property that each cone
singularity of $Z$ lies in the image of a vertex of $T$.  Isotope $h$
to send each edge of $T$ to its geodesic representative {\em
rel}-endpoints on $Z$ (if an edge $e$ terminates at a puncture, $h$
should send $e$ to a geodesic arc asymptotic to the corresponding cusp
of $Z$).  Then if $N$ is a hyperbolic 3-manifold and there is a
path-isometry $g \colon Z
\to N$ that is a local isometry on $Z - T$, then we call the pair
$(g,Z)$ a  
{\em simplicial hyperbolic surface in $N$} with
{\em associated triangulation $T$}.

Often the construction goes in the other direction: given a
triangulation $T$ of $\interior(S)$ with one vertex $v \in
\interior(S)$ and at least one edge $e$ so that $e \cup v$ forms a
closed loop $\gamma$, we can straighten any smooth, proper,
incompressible map $g'$ of $\interior(S)$ to $N$ to a simplicial
hyperbolic surface with associated triangulation $T$.  First we
straighten $g'$ so it maps $\gamma$ to its geodesic representative.
Then, straightening $g'$ on the edges of $T$ {\em rel}-endpoints
(possibly ideal endpoints) and then on faces of $T$, we obtain a map
$g \colon \interior(S) \to N$.  The pull-back metric from $N$
determines a singular hyperbolic surface $Z$ with a cone singularity
at $v$; since $v$ lies on a closed geodesic that is mapped to a closed
geodesic in $N$, the cone angle at $v$ is at least $2\pi$.  The result
is a simplicial hyperbolic surface $(g,Z)$ in $N$ with associated
triangulation $T$.  In this case we say that $(g,Z)$ is {\em adapted}
to $\gamma$.

Let $(f \colon S \to N) \in AH(S)$.  Then 
we denote by $\sh_k(N)$ the {\em marking preserving simplicial
hyperbolic surfaces in $N$ with at most $k$ cone-singularities}, namely,
simplicial hyperbolic surfaces $(g,Z)$, with $(h \colon \interior(S) \to Z) \in
\Sing_k(S)$, so that $g \compos h$ is homotopic to $f$.
If $\sigma$ is a simplex in $\calC(S)$ with vertices
$\alpha_1, \ldots, \alpha_p$, then we say a simplicial hyperbolic
surface $(g,Z) \in \sh_k(N)$ {\em realizes} $\sigma$ if $g$ maps each
$\alpha_i$ isometrically to its geodesic representative in $N$.  

Recall that a manifold $(f\colon S \to N) \in AH(S)$ has an {\em
accidental parabolic} if there is a non-peripheral element $\gamma \in
\pi_1(S)$ so that $f_*(\gamma)$ is a parabolic element of
$\Isom^+(\half^3)$.  Applying Hatcher's theorem
\cite{Hatcher:triangulations} allowing one to connect triangulations
by elementary moves, Canary proves the following (see \cite[\S
5, \S 6]{Canary:inj:radius}; compare \cite[\S 3]{Brock:length} and
\cite[\S 4]{Minsky:torus}):
\begin{theorem}[Canary]
Let $N \in AH(S)$ have no accidental parabolics, and let $(g_1,Z_1)$
and $(g_2,Z_2)$ lie in $\sh_1(N)$ where $(g_1,Z_1)$ is adapted to
$\alpha$ and $(g_2,Z_2)$
is adapted to $\beta$.  Then there is a continuous
family $(g_t \colon Z_t \to N) \subset \sh_2(N)$, $t \in [1,2]$.
\label{theorem:canary}
\end{theorem}

Using such an interpolation, we now give the proof of Lemma~\ref{lemma:count}.

\bold{Proof:} {\em (of Lemma~\ref{lemma:count}).}
Given the quasi-Fuchsian manifold $Q(X,Y)$, let $P_X$ and $P_Y$ denote
pants decompositions for which $X$ lies in 
$V(P_X)$ and $Y$ lies in $V(P_Y)$.  Then by a theorem of Bers
(\cite[Thm. 3]{Bers:bdry}, \cite[Prop. 6.4]{McMullen:iter}), we have 
$$\ell_{Q(X,Y)}(\alpha) < 2L$$ for each curve $\alpha$ in $P_X \cup
P_Y$.

Let $Z_X$ and $Z_Y$ denote simplicial hyperbolic surfaces realizing
$P_X$ and $P_Y$ in $Q(X,Y)$ and let $T_X$ and $T_Y$ denote their
associated triangulations.  
Let $Z^{\rm h}_X \in \Teich(S)$ be the hyperbolic surface
conformally equivalent to $Z_X$,
and let $Z^{\rm h}_Y$ be the hyperbolic surface conformally equivalent 
to $Z_Y$.  Finally, let $P_I$ and $P_T$ be pants
decompositions so that  
$Z^{\rm h}_X \in V(P_I)$ and $Z^{\rm h}_Y \in
V(P_T)$.

The next step will be to interpolate simplicial hyperbolic surfaces
between $Z_X$ and $Z_Y$ and estimate the minimum number of sets
$V(P_j)$ the corresponding conformally equivalent hyperbolic
representatives in $\Teich(S)$ intersect.  To show we have not
sacrificed too much distance in the pants graph, we prove the
following:
\begin{lem}
There is a constant $B' >0$ depending only on $S$ so that
$$\max\left\{d_{\bf P}(P_X,P_I), d_{\bf P}(P_Y,P_T)\right\}< B'.$$
\label{lemma:pants:comparison}
\end{lem}
\bold{Proof:} We have
$$\ell_{Z_X}(\gamma) <  2L$$ 
for each $\gamma \in 
P_X$.  
By a lemma of Ahlfors
\cite{Ahlfors:simplicial} we have
$$\ell_{Z^{\rm h}_X}(\beta) \ge \ell_{Z_X}(\beta)$$ for each
$\beta \in \eS$.  
Thus, we have 
$$\ell_{Z_X}(\alpha) < 2L$$ for each $\alpha \in P_X \cup
P_I$.  
Since curves in $P_X \cup P_I$ are
realized with bounded length on the simplicial hyperbolic surface $Z_X
\in \sh_k(Q(X,Y))$, 
we may apply 
\cite[Lem. 3.3]{Brock:length}
to find a $C'>0$ depending only on $L$ so that 
$$i(P_X, P_I) < C'.$$

Arguing similarly for $Z_Y$, we have
$$i(P_Y,P_T) < C'.$$
Applying the proof of Theorem~\ref{lemma:unit:length}, we have a
$B'>0$ depending only on $C'$ for which
$$\max\{d_{\bf P}(P_X,P_I), d_{\bf
P}(P_Y,P_T)\} <  B'.$$
\qed

To complete the proof of Lemma~\ref{lemma:count}, we seek
a continuous family 
$$(h_t \colon Z_t \to Q(X,Y)) \subset
\sh_k(Q(X,Y))$$ 
of simplicial hyperbolic surfaces in $Q(X,Y)$
interpolating between $(h_X,Z_X)$ and $(h_Y,Z_Y)$.  We first connect
$(h_X,Z_X)$ and 
$(h_Y,Z_Y)$ to simplicial hyperbolic surfaces $(h_X',Z_X')$ and
$(h_Y',Z_Y')$ adapted to single curves $\alpha \in P_X$ and $\beta
\in P_Y$ by continuous families, and then apply the interpolation
arguments of Canary.  

This is easily done by collapsing edges of $T_X$ that join distinct
vertices of $T_X$ down to a single vertex: if $e$ is such an edge
adjacent to a vertex $v$ on $\alpha$, then we may effect such a
collapsing by `dragging $h_X$ along $h_X(e)$' (see \cite[\S
5]{Canary:inj:radius}).  Precisely, if $v'$ is the other vertex in
$\bdry e$, we construct a homotopy of $(h_X,Z_X)$ to a new simplicial
hyperbolic surface by pulling the image $h_X(v')$ along the geodesic
segment $h_X(e)$ to $h_X(v)$ and pulling the edges and faces adjacent
to $v$ along with it while keeping the triangulation straight: the
image of each triangle is required to lift to the convex hull of its
vertices throughout the homotopy.

It is easy to check that under such a collapsing the cone angles at
the vertices remain at least $2 \pi$, and the number of vertices in
$T$ is reduced by 1.  We successively collapse edges joining
$v$ to different vertices until we are left with a triangulation with the
single vertex $v$, and a simplicial hyperbolic surface $(h_X',Z_X')
\in \sh_1(Q(X,Y))$ adapted to $\alpha$.

We perform analogous collapsings on the associated triangulation
for $(h_Y,Z_Y)$ to obtain the simplicial hyperbolic 
surface $(h_Y',Z_Y') \in \sh_1(Q(X,Y))$ adapted to $\beta$.  By
Theorem~\ref{theorem:canary}, we may
interpolate between $(h_X',Z_X')$ and 
$(h_Y',Z_Y')$ by a continuous family of simplicial hyperbolic
surfaces, so we have the desired continuous family 
$$(h_t \colon Z_t \to C(X,Y)) \subset
\sh_k(Q(X,Y)), \ \ \ t \in [0,1]$$ so that $(h_0,Z_0) = (h_X,Z_X)$ and
$(h_1,Z_1) = (h_Y,Z_Y)$. 

The singular hyperbolic structures $Z_t$ determine a
continuous path 
$$(h_t^{\rm h} \colon \interior(S) \to Z_t^{\rm h}) \subset \Teich(S),$$
where as before $Z^{\rm h}_t$ is the 
finite-area hyperbolic structure on $S$ in the same conformal class as
$Z_t$.  
Since $Z_t$ and $Z^{\rm h}_t$ represent metrics
on the same underlying surface $\interior(S)$, we have a natural continuous family of
1-Lipschitz mappings
$$\hat{h_t} \colon Z_t^{\rm h} \to C(X,Y)$$
of hyperbolic surfaces $Z_t^{\rm h} \in \Teich(S)$
into $Q(X,Y)$ so that for each $t$, $\hat{h_t}$ factors through the
simplicial hyperbolic surface $h_t \colon Z_t \to C(X,Y)$.

There are pants decompositions $P_1, \ldots, P_N$ (possibly with
repetition) that determine an open
cover $\{U_j\}_{j=0}^N$ of $[0,1]$ so that if $t \in 
U_j$ then $Z_t^{\rm h}$ lies in $V(P_j)$, and so that 
$U_j \cap U_{j+1} \not= \nullset$ for each $j = 0, \ldots, N-1$.
Applying Lemma~\ref{lemma:unit:length}, the sequence of pants decompositions
$$g = \{ P_I = P_0, \ldots, P_N = P_T \}$$
satisfies the hypotheses of Lemma~\ref{lemma:projection:argument} with 
$k = B.$  Applying Lemma~\ref{lemma:projection:argument}, we have a
$K_0$ so that $$d_{\bf P}(P_I,P_T) \le K_0 |\eS_g|.$$

Since the non-empty set $\eS_g$ is a subset of $\eS_{< L}$, and
Lemma~\ref{lemma:pants:comparison} guarantees 
$$d_{\bf P}(P_I,P_T) \ge d_{\bf P}(P_X,P_Y) - 2B',$$
we may combine constants to obtain a $K$ for which
$$d_{\bf P}(P_X,P_Y) \le K |\eS_{< L}|$$
proving the lemma.
\qed

Given a hyperbolic 3-manifold $M$, we let $\eG_{<L}(M)$ denote the set of
homotopy classes of closed geodesics in $M$ with length bounded above
by a constant $L >0$.  We use the 
contraction $\eG_{<L} = \eG_{<L}(M)$ when the manifold $M$ is understood.
The next lemma shows that the size $|\eG_{<L}|$ of $\eG_{<L}$ provides a lower
bound for the convex core volume of a hyperbolic 3-manifold in a general
context. 

\begin{lem}
Let $M$ be a geometrically finite hyperbolic 3-manifold with $\bdry M$
incompressible, and let $\vol(M)$ denote its convex core volume.  
Then there is a constant
$C_1>1$ depending only on $L$ and $C_2>0$ depending only on
$\chi(\bdry M)$ for which $$\frac{|\eG_{<L}|}{C_1} - C_2 < \vol(M).$$
\label{lemma:isotopy:classes}
\end{lem}

\bold{Proof:} Let $\epsilon>0$ be less than the minimum of $L/2$ and the
3-dimensional Margulis constant.  Let $\calV$ be any maximal set of
points in the $\epsilon$-thick part $\core(M)_{\ge \epsilon}$ of the
convex core $\core(M)$ for which points in $\calV$ are separated by a
distance at least $\epsilon/2$.  
Letting $B(x,R)$ denote the ball of radius $R$ about $x$ in $M$, 
it follows that $B(x,\epsilon/4)$ is embedded in
the $\epsilon$-neighborhood $\calN_\epsilon(\core(M)_{\ge \epsilon})$
of the $\epsilon$-thick part of $M$,   and 
$$B(x,\epsilon/4) \cap B(x',\epsilon/4) = \nullset$$ for $x \not= x'$
in $\calV$.
By maximality, however, we have $$\core(M)_{\ge \epsilon}
\subset \bigcup_{x \in
\calV}B(x,\epsilon/2).$$

Each isotopy class $\beta \in \eG_{<L}$ has a representative
$\beta^\star \subset \core(M)_{\ge \epsilon}$ with 
arclength less than $L$.  
Since the $\epsilon/2$-balls about points in $\calV$ cover
$\core(M)_{\ge \epsilon}$, each $\beta^\star$ intersects $B(x,\epsilon/2)$
for some $x\in \calV$.  Given $x \in \calV$, let $\calA_x \subset
\eG_{<L}$ denote the set 
$$\calA_x = \{ \beta \in \eG_{<L} \st \beta^\star \cap B(x,\epsilon/2) \not= 
\nullset \}.$$

Lifting to the universal cover so that $x$
lifts to the origin $0 \in \half^3$, the elements $\beta \in \calA_x$
determine pairwise disjoint translates of the ball
$B(0,\epsilon/2) \subset \half^3$  lying
within the ball $B(0,L + 2\epsilon) \in \half^3$.
It follows 
that the number of elements in each $\calA_x$ satisfies
$$|\calA_x| < \frac{\vol(B(0,L + 2\epsilon))}{\vol(B(0,\epsilon/2))}$$
which we set equal to $C_0$.

Since every $\beta \in \eG_{<L}$ lies in some $\calA_x$, we have
$$\frac{|\eG_{<L}|}{C_0} \le |\calV|.$$
Since the balls of radius $\epsilon/4$ about points $x \in \calV$ are
embedded and pairwise disjoint in $\calN_\epsilon(\core(M)_{\ge \epsilon}),$
we have the lower bound
\begin{equation*}
\frac{|\eG_{<L}|}{C_0}\cdot \vol(B(0,\epsilon/4)) 
\le \vol_\epsilon(M)
\end{equation*}
where $\vol_\epsilon(M) = \vol(\calN_\epsilon(\core(M))).$
There is a constant $K_\epsilon>0$ depending only on $\epsilon$ and
$M$ so that 
$$\vol_\epsilon(M) - \vol(M) < K_\epsilon$$
(see e.g. \cite[Lem. 8.2]{Canary:ends},
\cite[8.12.1]{Thurston:book:GTTM}). 
The lemma follows by setting $C_1 = C_0/\vol(B(0,\epsilon/4))$ and $C_2
= K_\epsilon$.
\qed

\bold{Proof:}{\em (of Theorem~\ref{theorem:volume:lower})}.
Since $\eS_{<L}$ is a subset of $\eG_{<L}(Q(X,Y))$, we may
combine Lemma~\ref{lemma:count}
with Lemma~\ref{lemma:isotopy:classes} to obtain
\begin{equation}
\frac{d_{\bf P}(P_X,P_Y)}{K \cdot C_1}
-C_2 \le
\vol(X,Y).
\end{equation}
Applying Theorem~\ref{theorem:pants:qi} we have
\begin{equation}
\frac{d_{\rm WP}(X,Y)}{A_1 \cdot K \cdot C_1}
- \frac{A_2}{K \cdot C_1} - C_2 \le
\vol(X,Y).
\end{equation}
Letting $$K_1 = A_1 \cdot K \cdot C_1 \ \ \
\text{and} \ \ \  K_2 = \frac{A_2}{K \cdot C_1} + C_2$$ 
the theorem follows.
\qed

\section{Bounding the core volume from above}
\label{section:upper}

Our goal in this section will be to prove the following theorem.
\begin{theorem}
Given $S$, there are constants $K_3$ and $K_4$ so that if $Q(X,Y) \in
QF(S)$ is a quasi-Fuchsian manifold and $P_X$ and $P_Y$ are pants
decompositions for which $X \in V(P_X)$ and $Y \in V(P_Y)$ then we
have $$\vol(X,Y) \le K_3 d_{\bf P}(P_X,P_Y) + K_4.$$
\label{theorem:volume:upper}
\end{theorem}
Given Theorem~\ref{theorem:volume:lower} and
Theorem~\ref{theorem:pants:qi}, Theorem~\ref{theorem:volume:upper}
represents the final step in the proof of Theorem~\ref{theorem:main}.

Let $G \subset {\bf P}(S)$ be a shortest path joining $P_X$ and $P_Y$
so that the length of $G$ is simply $d_{\bf P}(P_X,P_Y)$.  Let $${\rm
spin}(G) = \{ \alpha^* \st \alpha \in P, \ P \in G \}$$ denote the
geodesic representatives in $Q(X,Y)$ of elements of the pants
decompositions along $G$.  We call these geodesics the {\em spinning
geodesics} for $G$; they will serve to anchor various tetrahedra in
$Q(X,Y)$ at their vertices; we will then ``spin'' these tetrahedra by
pulling their vertices around the geodesics.

Our upper
bound for $\vol(X,Y)$ will come from a model manifold $N = S\times I$
comprised of blocks that are adapted to ${\rm spin}(G)$, together with
a piecewise $C^1$ surjective homotopy equivalence $f \colon N \to
\calN_\epsilon(C(X,Y))$ so that the image of each block under $f$ has
uniformly bounded volume.  The model will decompose into two parts.
\begin{enumerate}
\item {\bf The Caps:}  At each end of $N$ are {\em caps}, namely
products $S\times I$ on which $f$ restricts to homotopies of
simplicial hyperbolic surfaces $$h_X \colon Z_X \to C(X,Y) \ \
\text{and} \ \ h_Y \colon Z_Y \to C(X,Y)$$ realizing $P_X$ and $P_Y$
to the boundary components $X^\epsilon_h$ and $Y^\epsilon_h$ of the
$\epsilon$-neighborhood $\calN_\epsilon(C(X,Y))$ of the convex core.
\item {\bf The Triangulated Part:}  The caps sit at either end of the
{\em triangulated part} $N_\Delta$, a union of tetrahedra on which $f$
is simplicial: $f$ lifts to a map sending each simplex to the convex
hull of its vertices.  It follows that the image of each tetrahedron
$\Delta \in N_\Delta$ under $f$ has uniformly bounded volume.  We use
the geodesics $\alpha^*$, where $\alpha \in P \in G$, as a scaffolding
to build $N_\Delta$, a glueing of tetrahedra whose image interpolates
between the simplicial hyperbolic surfaces $Z_X$ and $Z_Y$.  After
``spinning'' $f$ sufficiently far about the spinning geodesics, all but 
a constant times $d_{\bf P}(P_X,P_Y)$ of the tetrahedra in $N_\Delta$
have images with small volume.
\end{enumerate}
These two arguments give the desired bound after collecting constants.

\bold{Remark:}  The above spinning trick is inspired by the ideal
simplicial maps of \cite{Thurston:hype1} which are in effect a limit
of the spinning process we perform here.  The result in our context of
passing to such a limit is an {\em ideal} triangulation of all but a
bounded volume portion of $C(X,Y)$, with a uniformly bounded number of
ideal tetrahedra necessary to accomplish each individual elementary
move (the small volume tetrahedra collapse to lower dimensional
ideal edges and faces).  We have chosen to work with finite
triangulations in the interest of demonstrating how the combinatorics
of ${\bf P}(S)$ may be used to produce triangulations of
3-manifolds an semi-algorithmic manner, independent of any geometric structure.

\subsection{Triangulations of surfaces} 
We specify a type of triangulation
of $S$ that is suited to a pants decomposition $P$.  By a {\em pair of 
pants} we will mean a connected component $\widehat{S}$ of $S - \calN(P)$, the
complement of the union of pairwise disjoint open annular
neighborhoods $\calN(P)$ of the curves in $P$ on $S$.
\begin{defn}
A {\em standard triangulation} $T(\widehat{S})$ for a pair of pants
$\widehat{S}$ is a triangulation with the following properties:
\begin{enumerate}
\item $T(\widehat{S})$ has two vertices on each boundary component.
\item $T(\widehat{S})$ has two disjoint {\em spanning triangles} with
no vertices in common, and a vertex on each component of $\bdry
\widehat{S}$.
\item The remaining 3 quadrilaterals are diagonally subdivided by an
arc that travels ``left to right'' with respect to the inward pointing 
normal to $\bdry \widehat{S}$ (see
Figure~\ref{figure:standard:triangulation}).
\end{enumerate}
\label{definition:standard:triangulation}
\end{defn}

\makefig{Standard triangulations suited to a pants decomposition.}{figure:standard:triangulation}{\psfig{file=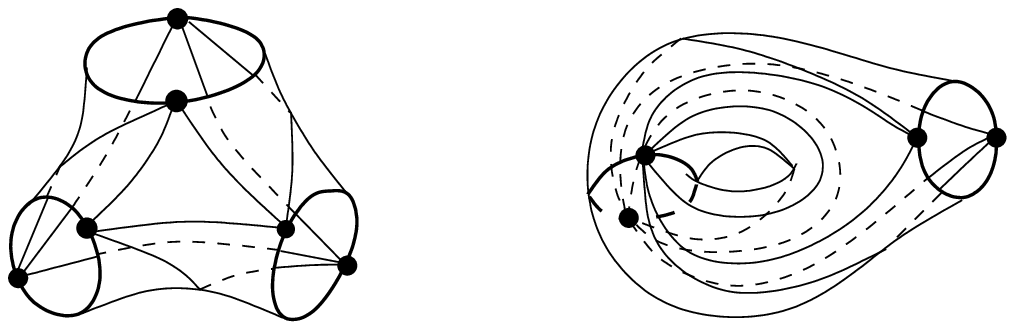,height=1.5in}}

We construct a standard triangulation suited to a pants decomposition $P
\in {\bf P}(S)$ by gluing together standard triangulations on pairs
of pants $\widehat{S}$ as follows.
\begin{defn}
Given a pants decomposition $P \in {\bf P}(S)$, a {\em standard
triangulation suited to $P$} is a triangulation $T$ of $S$ obtained as 
follows (see Figure~\ref{figure:standard:triangulation}):
\begin{enumerate}
\item $T$ has two vertices $p_\alpha$ and $\bar{p}_\alpha$ on each
component $\alpha$ of $P$, and two edges $e_\alpha$ and
$\bar{e}_\alpha$ in the complement $\alpha - p_\alpha \cup
\bar{p}_\alpha$.
\item If $S_0$ is a complementary open pair of pants in $S -P$, then
the restriction of $T$
to $S_0$ is the restriction to $S_0$ of a standard triangulation of
$S_0 \cup \bdry S_0$ (in the sense of Definition~\ref{definition:standard:triangulation}).
\item If two boundary
components $\alpha_1$ and $\alpha_2$ of $S_0$ are identified in $S$,
then the edge of each spanning triangle that runs from $\alpha_1$ to
$\alpha_2$ forms a closed loop.
\end{enumerate}
\label{defn:triangulation}
\end{defn}

\bold{Moves on triangulations.}
Given an elementary move on pants decompositions $(P,P')$,
i.e. $d_{\bf P}(P,P') = 1$, we 
now describe simple moves on triangulations that allow us to 
move from a standard triangulation suited to $P$ to a standard
triangulation suited to $P'$.  To distinguish moves on triangulations
from moves on pants decompositions, we refer to the latter as {\em pants
moves}.

To fix notation, given a pants move $(P,P')$ let $\alpha \in P$
and $\beta \in P'$ be the curves for which $i(\alpha,\beta) \not=
0$.  We call $\alpha$ and $\beta$ the curves {\em involved in the
pants move $(P,P')$}.  Let $\calN(P-\alpha)$ denote the union of
pairwise disjoint open annular neighborhoods about the curves in
$P-\alpha$, and let $S_\alpha$ denote the essential subsurface of
$S- \calN(P - \alpha)$ containing $\alpha$.  If $S_\alpha$ has genus 1 then
$(P,P')$ is called a {\em genus 1}  
pants move.  Likewise,  $S_\alpha$ has genus 0 then
$(P,P')$ is called a {\em genus 0} pants move.
We say the pants move $(P,P')$ occurs on $S_\alpha$.

A standard triangulation $T$ suited to $P$ naturally identifies
candidate elementary moves for each $\alpha \in P$: there is a natural
choice of isotopy class of simple closed curves $\beta \subset
S_\alpha$ for which $i(\alpha,\beta) =1$ or $i(\alpha,\beta) = 2$
depending on whether $S_\alpha$ has genus 1 or genus 0.  If $S_\alpha$
has genus 1, then each spanning triangle for $T$ in $S_\alpha$ has one
edge with its endpoints identified.  These edges are in the same
isotopy class which we call $\beta(\alpha,T)$.  Likewise, if
$S_\alpha$ has genus $0$, then removing the edges of $T$ that do not
have endpoints lying on $\alpha$ produces two hexagons in the
complement of the remaining edges.  Concatenating edges in these
hexagons joining the two vertices in each that lie on $\alpha$ we
obtain an isotopy class of simple closed curves, which we again call
$\beta(\alpha,T)$.

For each $\alpha \in P$, the pants decomposition $P' = (P-\alpha) \cup
\beta(\alpha,T)$ satisfies $d_{\bf P}(P,P') = 1.$

\medskip

There are three basic types of moves on these triangulations:

\ital{{\bf MVI.} The Dehn twist move.}  One standard move on triangulations we
will use effects a Dehn twist of a standard triangulation suited to
$P$ about a curve $\alpha \in P$.  Given a triangulation $T$ suited to
$P$, let $TW_\alpha(T)$ denote the standard triangulation suited to
$P$ obtained by shifting each edge with a vertex on $\alpha$ to the
right along $\alpha$ until it hits the next vertex.  Then
$TW_\alpha(T)$ is isotopic to the image of $T$ under a right-$\alpha$
Dehn twist.  We define $TW_\alpha^{-1}(T)$ similarly, by shifting
edges to the left rather than to the right.

Note that for any triangulation $T$
suited to $P$ and any $\alpha \in P$ we have
$$\beta(\alpha,TW_\alpha(T)) = \tau_\alpha(\beta(\alpha,T))$$
where $\tau_\alpha$ is a right $\alpha$-Dehn twist.

The other elementary moves on triangulations will be specific to a
given type of pants move.   Given a triangulation $T$ suited to
$P$ and $\alpha \in P$, we describe blow-down and blow-up moves that
allow us to pass from a triangulation suited to $P$ to a triangulation 
suited to the pants decomposition $(P - \alpha) \cup \beta(\alpha,T)$.

\makefig{The genus 1 blow-down and blow-up moves.}{figure:grafting}{\psfig{file=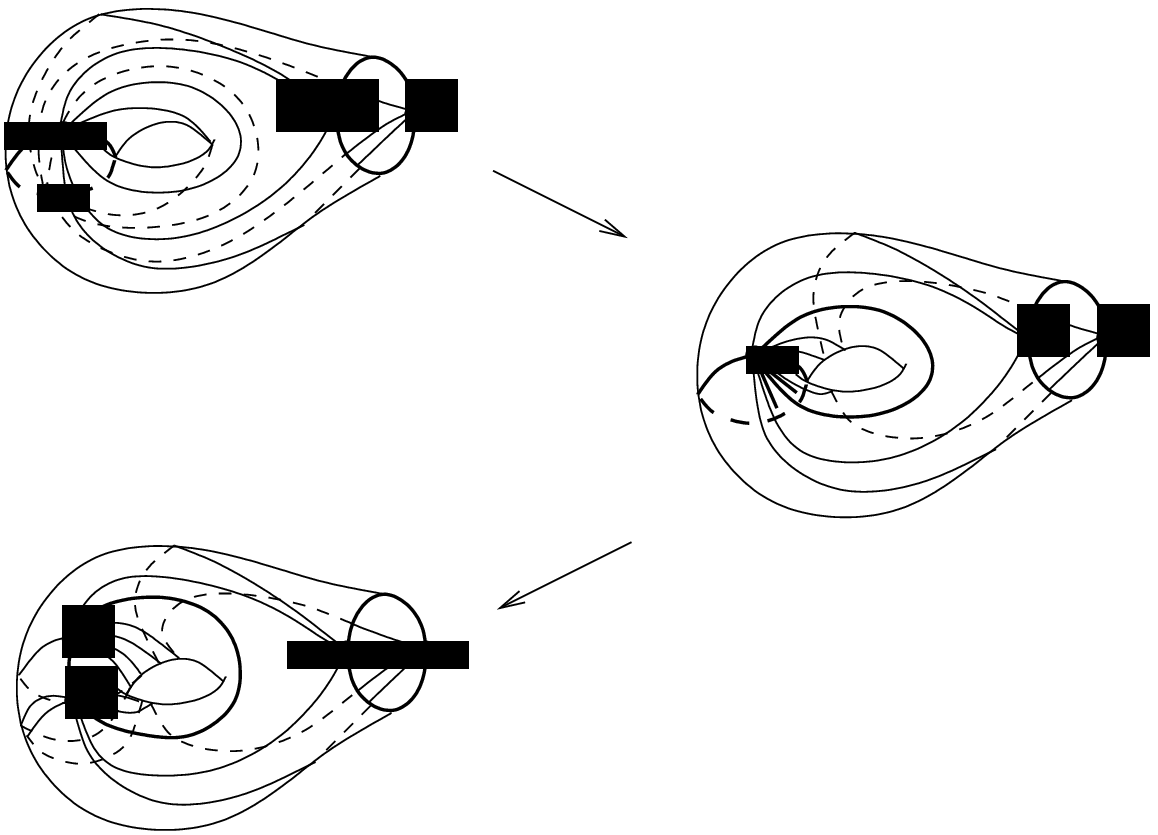,height=2.4in}}
\ital{{\bf MVII.} Genus 1 moves.}  
Given a standard triangulation $T$ suited to $P$
the edges of $T$ that close to form loops in the isotopy class of
$\beta(\alpha,T)$ bound an
annulus on $S$.  We call this annulus $A_\alpha$ the {\em $\alpha$
compressing annulus} for $T$.  The {\em genus 1 $\alpha$-blow down}
$BD_\alpha^1(T)$ of $T$ is the triangulation of $S$ obtained by
collapsing the arc $e$ of $\alpha 
\cap A_\alpha$ to a point and collapsing the two triangles in $T$
containing $e$ to a single edge.

Let $\beta = \beta(\alpha,T)$.  Then given $T' = BD_\alpha^1(T)$
the {\em genus 1 $\beta$-blow-up} $BU_\beta^1(T')$ of $T'$, is obtained
by grafting a compressing annulus $A_\beta$ in along the curve $\alpha$
to obtain a standard triangulation suited to $(P - \alpha) \cup
\beta$ as in Figure~\ref{figure:grafting}.

\ital{{\bf MVIII.} Genus 0 moves.}  Let $T$ be a
standard triangulation suited to $P$.  Let $e_\alpha$ and
$\bar{e}_\alpha$ denote the two edges of $T$ that constitute the curve
$\alpha$.  Call these curves the $\alpha$-edges of $T$.  As described
above there are two hexagons $H$ and $\bar{H}$ obtained by removing
all edges of $T$ that do not have endpoints on $\alpha$.
The {\em genus 0 $\alpha$-blow-down} $BD_\alpha^0(T) = T'$
(Figure~\ref{figure:g0bd}) is obtained by
performing 3 ``diagonal switches'' on each hexagon
(Figure~\ref{figure:hexagons}) 
\makefig{Diagonal switches on hexagons.}{figure:hexagons}{\psfig{file=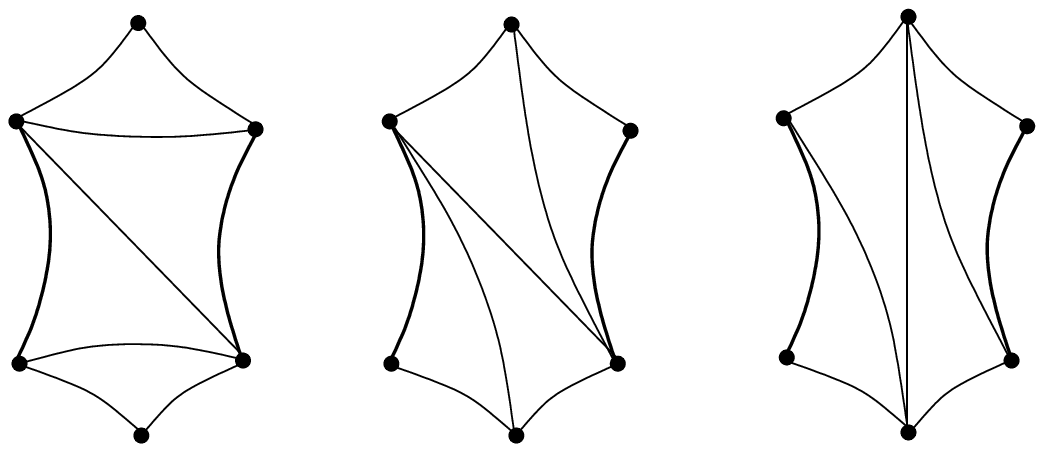,height=1.8in}} 
to yield a new triangulation with edges 
$e_\beta$
and $\bar{e}_\beta$ as edges whose concatenation gives the curve
$\beta = \beta(\alpha,T)$. 

Likewise, the {\em genus 0 $\beta$-blow-up} $BU^1_\beta(T')$ is
obtained by modifying the hexagons containing the $\alpha$-edges by
the inverses of the 3-diagonal switches.  Note that $BU_{\beta}^0
\circ BD_\alpha^0(T)$ is a standard triangulation 
suited to the pants decomposition $(P - \alpha) \cup \beta$.

\makefig{The genus 0 blow-down and blow-up
moves.}{figure:g0bd}{\psfig{file=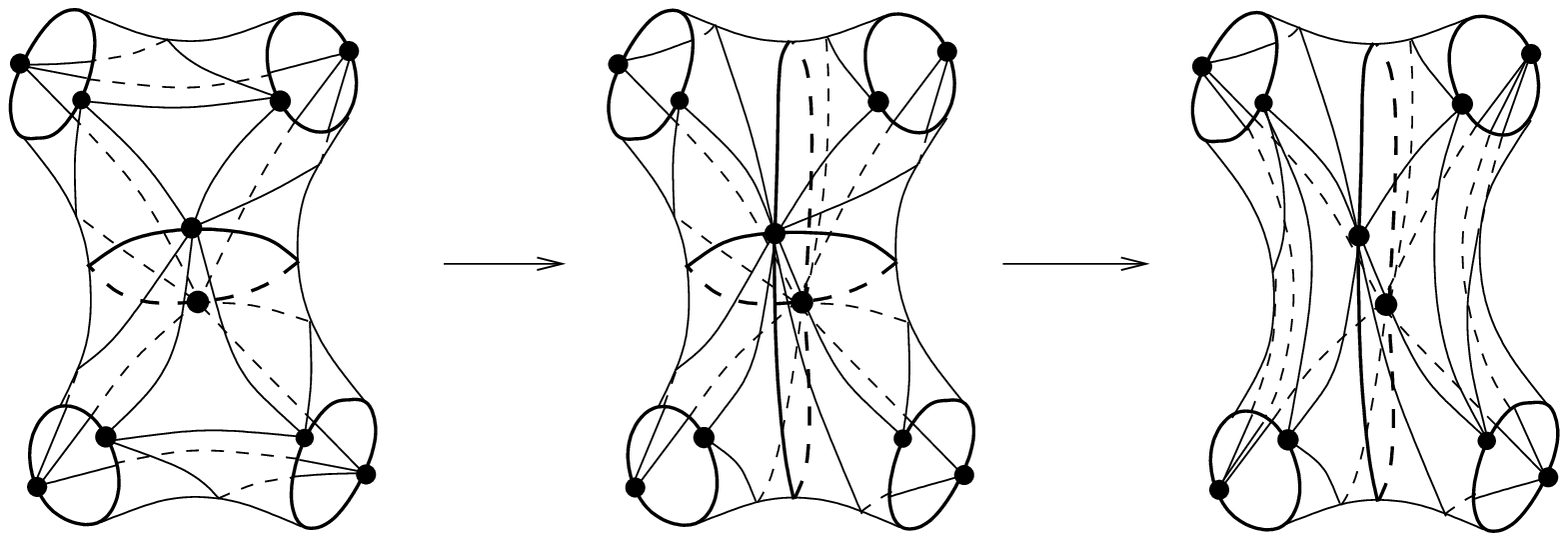,height=1.6in}} 

We summarize properties of these moves as a lemma.
\begin{lem}
Let $(P,P')$ be a pants move involving $\alpha \in P$ and $\alpha' 
\in P'$.  Let $T$ be a standard triangulation suited to $P$.  Then
there is an $n \in \zed$ so that 
$$BU_{\alpha'} \circ BD_\alpha \circ TW_\alpha^n(T)$$
is a standard triangulation suited to $P'$.
\label{lemma:triangulation:moves}
\end{lem}

\subsection{Realizing moves by blocks}  
As before let $Q(X,Y)$ be a quasi-Fuchsian 
manifold, and let $P_X$ and $P_Y$ be pants decompositions for which $X
\in V(P_X)$ and $Y \in V(P_Y)$.  Let $G \subset {\bf P}(S)$ be a
geodesic joining $P_X$ and $P_Y$.  Recall we denote by
${\rm spin}(G) = \{ \alpha^* \st \alpha \in P,\ P \in G\}$
the {\em spinning geodesics} associated to $G$.

Equip each spinning geodesic $\alpha^*$ with a pair of antipodal
vertices $p_\alpha$ and $\bar{p}_\alpha$: i.e. points on
$\alpha^*$ so that the distance from $p_\alpha$ to
$\bar{p}_\alpha$ along $\alpha^*$ is maximal.  For reference, we
equip each $\alpha$ and thence each $\alpha^*$ with an orientation.

Let $P_i$, $i=0, \ldots, m$, denote the pants decompositions along the
geodesic $G$, so that $P_0 = P_X$ and $P_m = P_Y$.
Making an initial choice of standard triangulation $T_X$ suited to $P_X$,
Lemma~\ref{lemma:triangulation:moves} provides a sequence of moves on
triangulations allowing us to process from the triangulation $T_X$ to a
standard triangulation $T_Y$ suited to $Y$ via standard triangulations
$T_i$ suited to $P_i$.  We begin with a model manifold 
$$N_0 \cong S \times I$$ and triangulate $\bdry^+ N_0 = S \times
\{0\}$ by $T_X$.  Our 
aim is to build models $N_i \cong S\times I$ by gluing triangulated
I-bundles to $\bdry^+ N_{i-1}$ so that the resulting triangulation on
$\bdry^+ N_i = T_i$.  We do this by building a triangulated subsurface
block corresponding to each elementary move and successively gluing
the blocks to $\bdry^+ N_{i-1}$.

\begin{defn}  Given a curve $\alpha$ in a pants decomposition $P$, 
we define the {\em subsurface block} by the quotient
$$B_\alpha = S_\alpha \times [0,1]/(x,t) \sim (x,0) \ \ \ \text{for} \ 
\ \ x \in \bdry S_\alpha, \ \ t \in [0,1].$$
We denote the {\em upper} and {\em lower} boundary of $B_\alpha$ by
$$\bdry^+B_\alpha = S\times\{0\} \ \ \ \text{and} \ \ \
\bdry^-B_\alpha = S\times\{1\}.$$
\end{defn}

We now describe {\em block triangulations} associated to each
elementary move on triangulations.  We will say a triangulation
$\calT$ of a subsurface block $B_\alpha$ {\em realizes} an elementary
move $T \to \calM(T)$ if we have
$$\calT \cap \bdry^- B_\alpha = T  \ \ \ \text{and} \ \ \ \calT \cap
\bdry^+ B_\alpha = \calM(T).$$

\bold{Block triangulations.}  
Let $T$ be a
standard triangulation on $S_\alpha$ suited to $\alpha$.  Then the
{\em standard block triangulation} $\calT_\alpha$ is obtained from 
$$T \times [0,1] / \sim$$ in the following way.  Initially, 
$T \times [0,1]/\sim$ is a cell decomposition of $B_\alpha$.  
For any edge $e$ of $T$ with $\bdry e \cap \bdry S_\alpha = \nullset$,
$e \times [0,1]$  is a quadrilateral to which we add a diagonal
depending on the genus of $S_\alpha$.
\begin{itemize}
\item When $S_\alpha$ has genus $0$, the $\alpha$-edges
$e_\alpha$ and $\bar{e}_\alpha$ determine quadrilaterals 
$e_\alpha \times [0,1]$ and $\bar{e}_\alpha \times [0,1]$ in $\calT$.
We triangulate these quadrilaterals with two new edges that run in the same
direction along the annulus 
$$e_\alpha \times 
[0,1] \cup \bar{e}_\alpha \times [0,1]$$ and cone off the new edges
down to the vertices opposite the quadrilaterals (see
Figure~\ref{figure:std:block}). 
\makefig{Extending $T\times [0,1]/\sim$ to a
triangulation.}{figure:std:block}{\psfig{file=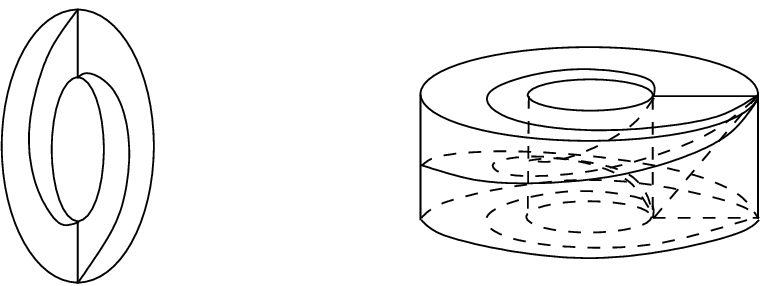,height=1.4in}} 
\item When $S_\alpha$ has genus $1$ in addition to the
two edges $e_\alpha$ and $\bar{e}_\alpha$ that concatenated give
$\alpha$ there are $3$ other edges 
$e_1$, $e_2$ and $e_3$ that triangulate the $\alpha$-compressing
annulus $A_\alpha$ for which $\bdry e_j \cap \bdry S_\alpha = \nullset$, $j =
1,2,3$.  We triangulate the annulus 
$$e_\alpha \times 
[0,1] \cup \bar{e}_\alpha \times [0,1]$$ as before, and extend this to 
a triangulation of $A_\alpha \times [0,1]$ with no new vertices as in
Figure~\ref{figure:std:block}.
\end{itemize}

From the standard block triangulation, we build four types of blocks:

\ital{{\bf BLI.}  The Dehn twist block.}
Given $S_\alpha$, the block triangulation $\calT({{TW}_\alpha})$ realizing the
move $TW_\alpha$ is obtained form the standard block triangulation
$\calT_\alpha$ as follows.  Consider the annulus $A  = \alpha \times [0,1]$ 
in $B_\alpha$ with the triangulation $T_A$ on $A$ induced  by
$\calT_\alpha$.
The reference orientation for $\alpha$ locally determines a {\em left} 
and {\em right} side of $\alpha$ in $S_\alpha$, and hence a left and
right side of $A$ in $B_\alpha$.  Cut $B_\alpha$ along $A$ to obtain
two annuli $A_L$ and $A_R$ that bound the {\em local} left and right
side of $B_\alpha - A$.  Re-glue the $\alpha \times \{0\}$ boundary
components of $A_L$ and $A_R$ by the identity, and re-glue the $\alpha 
\times \{1\}$ boundary components of $A_L$ and $A_R$ shifted by a
right Dehn twist.  

The triangulations induced by $T_A$ on $A_L$ and $A_R$ determine a
triangulation of the torus $A_L \cup A_R$ after re-gluing.  This
triangulation naturally extends to a triangulation of a solid torus
$V$ with boundary $A_L \cup A_R$ by filling in tetrahedra (the
triangulations on $A_L$ and $A_R$ differ by two pairs of diagonal
switches.  See Figure~\ref{figure:diagonal:switches}). After filling
in by $V$, the result is a standard block $B_\alpha$ with
triangulation $\calT_\alpha^{\rm tw+}$ realizing the move $TW_\alpha$.
We call the standard block $B_\alpha$ equipped with the triangulation
$\calT_\alpha^{\rm tw+}$ the {\em right Dehn twist block}
$B_\alpha^{\rm tw+}$.
\makefig{Triangulating the Dehn twist
block.\\ {\small The dotted edges indicate diagonal switches for $T_A$
and $tw(T_A)$ that produce a common triangulation of the annulus,
which results in the pictured triangulation of a solid torus
$V$.}}{figure:diagonal:switches}{\psfig{file=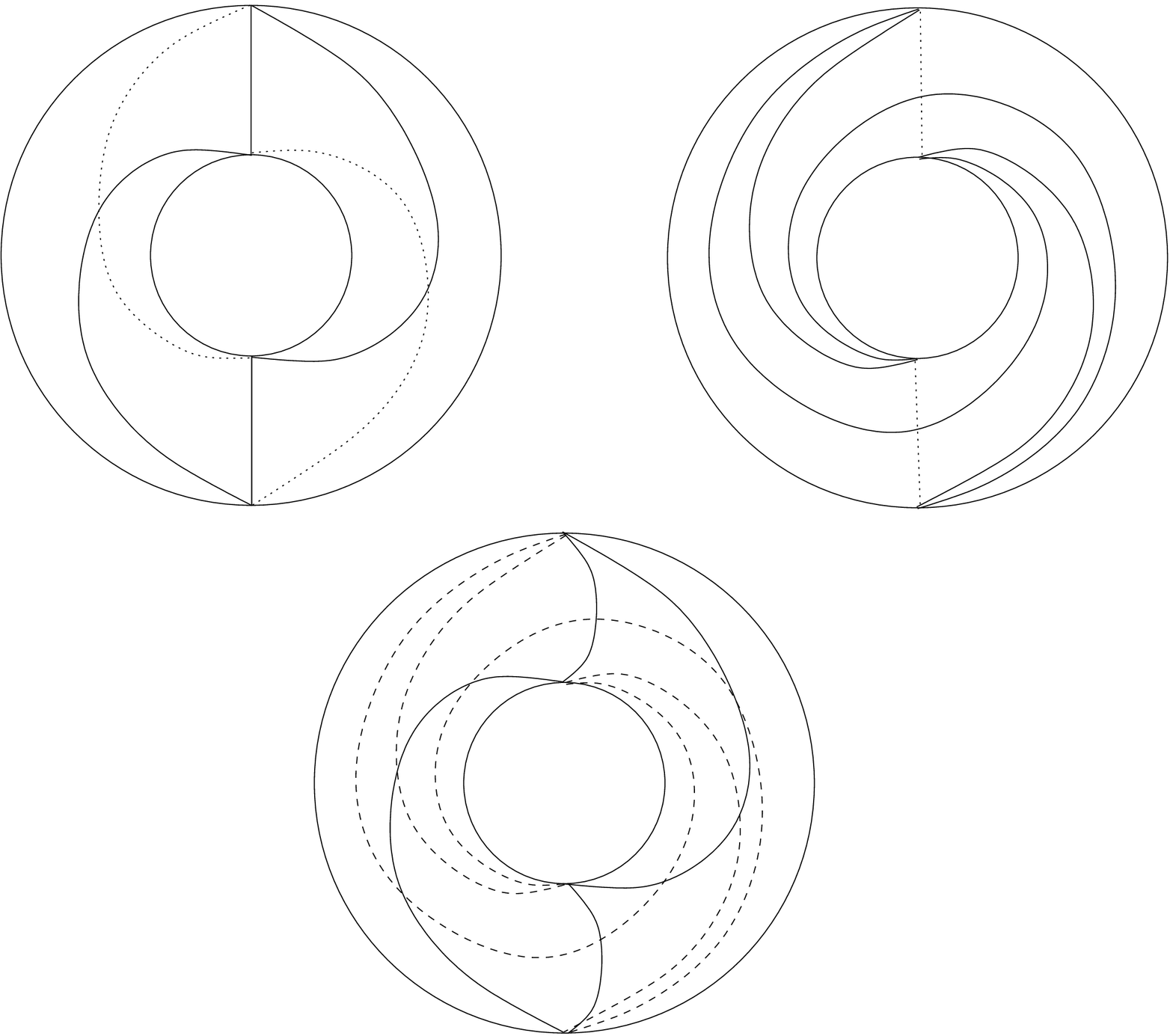,height=2.5in}}
The {\em left Dehn twist block} $B_\alpha^{\rm tw-}$ is obtained
analogously, by re-gluing with a left Dehn twist instead of a right
Dehn twist.  The left Dehn twist block $B_\alpha^{\rm tw-}$ carries
the triangulation $\calT_\alpha^{\rm tw-}$ and realizes the move
$TW_\alpha^{-1}$.

Blocks realizing the other moves are triangulated as follows.

\ital{{\bf BLII.} Blow-up and blow-down blocks.}  
We modify the standard block triangulation to obtain triangulated
blocks that realize blow-up and blow-down moves as follows.
To realize a genus 1
$\alpha$-blow-down by a block triangulation,
we modify the 
standard block triangulation $\calT_\alpha$ on $B_\alpha$ by
collapsing the compression annulus on $\bdry^+ B_\alpha$ as in the
description of the move $BD_\alpha^1.$  The only difference is that
here in addition to collapsing triangles to edges we also collapse
tetrahedra to triangles.  

Similarly, if $\beta = \beta(\alpha,T)$, we obtain a block
triangulation realizing the genus 1 
$\beta$-blow-up $BU_\beta^1$ by
collapsing the $\beta$ compression annulus on $\bdry^- B_\beta$.

The genus 0 $\alpha$-blow-down (and blow-up)
moves come from diagonal switches on hexagons.  We realize the move
$BD_\alpha^0$ by gluing tetrahedra realizing each of the diagonal
switches to $\bdry^+ B_\alpha$.
Likewise, the block triangulation realizing the genus 0
$\beta$-blow-up can be obtained by gluing tetrahedra in this way to
the $\beta$-hexagons on $\bdry^- B_\alpha$.  

We denote the $\alpha$-blow-up block by $B_\alpha^{\rm bu}$ and the
$\alpha$-blow-down block by $B_\alpha^{\rm bd}$.

\ital{{\bf BLIII.} Straightening blocks.}
There are two other types of block triangulations we will need that
realize the identity move on a blown-down triangulation.  We call these
triangulated blocks {\em straightening blocks}.  Given the
triangulation $BD_\alpha^1(T)$ or $BD_\alpha^0(T)$ for a standard
triangulation  $T$ suited to a pants decomposition $P$ containing $\alpha$, the
genus $g$ straightening triangulation, $g = 0,1$ is obtained by
completing the cell decomposition 
$$(BD_\alpha^g(T)\cap S_\alpha) \times [0,1]/\sim$$ of $B_\alpha$, and
extending this decomposition to a triangulation of $B_\alpha$ with no
new vertices. It is easy to check that this can be done.  We denote
the straightening blocks by $B_{\alpha,\beta}^{\rm st}$ where $\beta = 
\beta(\alpha,T)$.

\subsection{Mapping in blocks and building the model}
We now use our block triangulations to build a {\em model manifold} $N 
\cong S \times I$.  We build $N$ in stages corresponding to standard
triangulations $T_j$ suited to pants decompositions $P_j$ that intervene
between $P_X$ and $P_Y$ respectively.  At each stage
there is a model $N_j$ also homeomorphic to $S \times I$ so that the
top boundary component $\bdry^+ N_j = S \times \{1\}$ is triangulated
by the $j$th triangulation $T_j$ in the sequence of triangulations.
Each $N_j$ will be obtained from $N_{j-1}$ by 
attaching a triangulated block to $\bdry^+N_{j-1}$ that realizes the
elementary move needed to move from $T_{j-1}$ to $T_j$.

The model $N$ will come equipped with a map $f \colon N \to C(X,Y)$
that is {\em simplicial} on each block:  
\begin{defn}
An incompressible mapping $f \colon N \to M$ from a triangulated
3-manifold to a hyperbolic $3$-manifold is {\em simplicial} if the
lift $\widetilde{f} \colon \widetilde{N} \to \half^3$ sends each
$k$-simplex to the convex hull of its vertices.
\end{defn}

The following theorem describes the properties of our model and its
mapping to the quasi-Fuchsian manifold $Q(X,Y)$.  For simplicity we
assume for the remainder of this section that $S$ is closed, and
detail the necessary modifications to the argument at the end.

\begin{theorem}  Given $Q(X,Y)$, there is a model manifold $N \cong S
\times I$, 
equipped with a surjective homotopy equivalence
$$f \colon N \to \calN_\epsilon(C(X,Y))$$
to the $\epsilon$-neighborhood of the convex core $C(X,Y)$ with the following
properties: 
\begin{enumerate}
\item $N$ is the union $$N = {\rm cap}_X \cup N_\Delta \cup {\rm
cap}_Y$$
of the {\em caps} ${\rm cap}_X \cong S \times I$ and ${\rm cap}_Y
\cong S\times I$ and the {\em
triangulated part} $N_\Delta$, a union of 
blocks of the above type glued top boundary to bottom boundary, all
but $3 d_{\bf P}(P_X,P_Y)$ of which are {\em 
Dehn twist} blocks.
\item The map $f$ is piecewise $C^1$ and is simplicial on $N_\Delta$.
\item For each tetrahedron $\Delta$ in $N_\Delta$ that lies in a
Dehn twist block, $f$ maps some edge of $\Delta$ to a spinning geodesic.
\item The restriction $f \vert_{{\rm cap}_X}$ is a homotopy from
$\bdry^-\calN_\epsilon(C(X,Y))$ to a simplicial hyperbolic surface
realizing $P_X$, and the restriction $f\vert_{{\rm cap}_Y}$ is a
homotopy from a simplicial hyperbolic surface realizing $P_Y$ to
$\bdry^+ \calN_\epsilon(C(X,Y))$.
\end{enumerate}
\label{theorem:model}
\end{theorem}

\bold{Proof:} To motivate the construction of $N$ we build the map $f$
in stages as well. 

\ital{Mapping in ${\rm cap}_X$.}
Let $T_X$ and $T_Y$ be standard triangulations
suited to $P_X$ and $P_Y$.  Let
$$f_X \colon S\times [0,1] \to C(X,Y)$$ 
determine a homotopy of the convex core boundary 
$g_X \colon X_h \to \bdry^- C(X,Y)$ to a simplicial hyperbolic surface
$h_X \colon Z_X \to C(X,Y)$ with associated triangulation $T_X$ so
that if $v$ and $\bar{v}$ are the vertices of $T_X$ on $\alpha \in
P_X$, $f_X(v) = p_\alpha$ and $f_X(\bar{v}) = \bar{p}_\alpha$.
Notice that this implies that $h_X$ realizes each curve $\alpha \in P_X$.
Let $N_0 = S \times [0,1]$,  denote the domain for $f_X$; we will
refer to $N_0 = {\rm cap}_X$ as the $X$-cap of $N$. The top boundary
component $\bdry^+ N_0$ carries the triangulation $T_X$. 

Working inductively, we assume given a model $N_j$ at the $j$th stage:
i.e. 
\begin{enumerate}
\item the model $N_j \cong S \times I$, consists of the $X$-cap and a
triangulated part so that the upper boundary $\bdry^+ N_j$ is triangulated 
by a standard triangulation $T_j$ suited to $P_j$,
\item $N_j$ comes equipped with a map $f_j \colon N_j \to C(X,Y)$ that 
is simplicial on the triangulated  part of $N_j$,
\item $f_j \vert_{\bdry^+ N_j}$ factors through a simplicial 
hyperbolic surface $h_j \colon Z_j \to C(X,Y)$ with associated
triangulation $T_j$ so that $h_j$ sends vertices $v_\gamma$ and
$\bar{v}_\gamma$ on $\gamma \in P_j$ to $p_\gamma$ and
$\bar{p}_\gamma$ on $\gamma^*$.
\end{enumerate}
Let $\alpha \in P_j$ and $\alpha' \in P_{j+1}$ be the curves involved
in the genus $g$ elementary move $(P_j,P_{j+1})$, and let $n$ be the integer
guaranteed by 
Lemma~\ref{lemma:triangulation:moves} for which 
$$BU^g_{\alpha'} \circ BD_\alpha^g \circ TW_\alpha^{n}(T_j) =
T_{j+1}$$
is a standard triangulation suited to $P_{j+1}$.

We now specify how to add triangulated blocks to $N_j$ and extend
$f_j$ simplicially over each additional block to obtain the next stage 
of the model $N_{j+1}$.

\ital{Mapping in Dehn twist blocks.}  
If $n$ is positive, we attach $n$ right Dehn twist blocks to the model and
extend $f_j$ over them in sequence, while if $n$ is negative we do
likewise with left Dehn twist blocks.  

We assume $n$ is positive; the negative case is identical with left
Dehn twist blocks replacing right Dehn twist blocks.  We attach a
right Dehn twist block $B^{\rm tw+}_\alpha$ to $\bdry^+N_j$ so that the
triangulation on $\bdry^- B^{\rm tw+}_\alpha$ agrees with $T_j \cap
S_\alpha$ to obtain a new model $N_{j,1}$.  We extend $f_j$ over
$B_\alpha^{\rm tw+}$ to obtain a map $f_{j,1}$ as follows.  Recalling
that $B_\alpha^{\rm tw+}$ has the form $S_\alpha \times [0,1]/\sim$, we
set $f_j(x,t) = f_j(x,0)$.  We then straighten $f_j$ on $\alpha \times
\{1\}$ to its geodesic representative so that the vertices
$v$ and $\bar{v}$ on $\alpha \times \{1\}$ to $p_\alpha$ and
$\bar{p}_\alpha$, and finally we straighten $f_{j,1}$ by a homotopy
to make it simplicial on $B_\alpha^{\rm tw+}$.   
We note that every tetrahedron in the Dehn twist block  
has an edge that maps to a geodesic arc of
$\alpha^*$, verifying part (3). 

\smallskip
The map $f_{j,1}$ factors through a simplicial hyperbolic surface
still realizing $P_j$ with associated triangulation $TW_\alpha(T_j)$
which we denote by $T_{j,1}$.  Repeating this procedure to add $n$
Dehn twist blocks we arrive at a 
model $N_{j,n}$ equipped with a map $f_{j,n} \colon N_{j,n} \to
C(X,Y)$ so that 
\begin{enumerate}
\item $\bdry^+ N_{j,n}$ carries the triangulation 
$$T_{j,n} = TW_\alpha^n(T_j),$$
\item $f_{j,n}$ factors through a simplicial hyperbolic surface
realizing $P_j$ with associated triangulation $T_{j,n}$, and 
\item the vertices of $T_{j,n}$ map to $p_\gamma$ and
$\bar{p}_\gamma$ on $\gamma^*$ for each $\gamma \in P_j$.
\end{enumerate}

\ital{Mapping in blow-down blocks.}  Our discussion of how to attach a
blow-down 
block to $N_{j,n}$ and how to extend $f_{j,n}$ over this block breaks
into cases as usual.

\ital{Genus 0.}  The genus $0$ blow-down block $B_\alpha^{\rm bd}$ is
attached to $\bdry^+ N_{j,n}$ along $S_\alpha$ so that the triangulations
agree as before.  This gives a new model $N_{j,{\rm bd}}$.  We extend
$f_{j,n}$ over $B_\alpha^{\rm bd}$ to give a map $f_{j,{\rm bd}}
\colon N_{j,{\rm bd}} \to C(X,Y)$ by first mapping $B_\alpha^{\rm bd}$ 
to $\bdry^+ N_{j,n}$, as with the Dehn twist block, and then
straightening $f_{j,{\rm bd}}$ to a simplicial map.  

\ital{Genus 1.}  In the genus 1 case is the same, except that as there 
is only one vertex $v$ on $\alpha$ in $\bdry^+ B_\alpha^{\rm bd}$, we
simply send $v$ to $p_\alpha$ and straighten $f_{j,{\rm bd}}$ to a
simplicial map as before.

In each case the resulting map $f_{j,{\rm bd}} \vert_{\bdry^+
N_{j,{\rm bd}}}$ factors through a simplicial hyperbolic surface that
realizes $P_j$ and has associated triangulation
$BD_\alpha^0(T_{j,n})$ which we denote by $T_{j,{\rm bd}}$.

\ital{Mapping in straightening blocks.}  Straightening blocks allow us to pass
from simplicial hyperbolic surfaces realizing $P_j$ to simplicial
hyperbolic surfaces realizing $P_{j+1}$.  
We attach the straightening block $B_{\alpha,\alpha'}^{\rm st}$ to
$N_{j,{\rm bd}}$ to obtain a model $N_{j,{\rm st}}$.  We extend
$f_{j,{\rm bd}}$ over $B_{\alpha,\alpha'}^{\rm st}$ to a map
$f_{j,{\rm st}}$ by defining $f_{j,{\rm st}}(x,t) = f_{j,{\rm
st}}(x,0)$, and then straightening $f_{j,{\rm st}}$ on $\alpha' \times 
\{1\} \subset \bdry^+ B_{\alpha,\alpha'}^{\rm st}$.  We send the $\alpha'$
vertex  to $p_{\alpha'}$ or vertices to $p_{\alpha'}$ and
$\bar{p}_{\alpha'}$ and straighten the map to a simplicial map on
$B_{\alpha, \alpha'}^{\rm st}$.  Now
$$f_{j,{\rm st}} \vert_{\bdry^+ N_{j,{\rm st}}}$$
factors through a simplicial hyperbolic surface realizing $P_{j+1}$
with associated triangulation $T_{j,{\rm bd}}$ once again.

\ital{Mapping in blow-up blocks.}  This procedure is essentially the inverse of
the attaching and mapping in the blow-down blocks.

\ital{Genus 0.}  Attach $B_{\alpha'}^{\rm bu}$ to $\bdry^+ N_{j,{\rm
st}}$ along $S_{\alpha'}$ (which equals $S_\alpha$), and extend
$f_{j,{\rm st}}$ to $f_{j,{\rm bu}}$ over $B_{\alpha'}^{\rm bu}$ by
setting $f_{j,{\rm bu}}(x,t) = f_{j,{\rm bu}}(x,0)$ and then
straightening $f_{j,{\rm bu}}$ {\em rel}-$\alpha' \times \{1\}$ to a
simplicial map.

\ital{Genus 1.}  
Set $f_{j,{\rm bu}}(x,t) = f_{j,{\rm bu}}(x,0)$ and then homotope
$f_{j,{\rm bu}}$ so that it still sends $\alpha' \subset \bdry^+
B_{\alpha'}^{\rm bu}$ to its geodesic representative but also sends
the $\alpha'$ vertices $v_{\alpha'}$ and $\bar{v}_{\alpha'}$ to
$p_{\alpha'}$ and $\bar{p}_{\alpha'}$.

\medskip

The map $f_{j,{\rm bu}} \vert_{\bdry^+ N_{j,{\rm bu}}}$ factors through 
a simplicial hyperbolic surface realizing $P'$ with associated
triangulation 
$$T_{j+1} = BU_{\alpha'} \circ BD_\alpha \circ TW_\alpha^n(T_j)$$
which is a standard triangulation suited to $P_{j+1}$.  This completes
the inductive step.

\ital{Mapping in ${\rm cap}_Y$.}  Let $|G| = d_{\bf
P}(P_X,P_Y)$ denote the length of a geodesic $G \subset {\bf P}(S)$.
Then the above inductive procedure results finally in a map $$f_{|G|}
\colon N_{|G|} \to C(X,Y)$$ so that the restriction $f_{|G|,0}
\vert_{\bdry^+ N_{|G|}}$ factors through a simplicial hyperbolic
surface realizing $P_Y$ with associated triangulation $T_Y$ (a
standard triangulation suited to $P_Y$), and so that the vertices
$v_\gamma$ and $\bar{v}_\gamma$ map to the vertices $p_\gamma$ and
$\bar{p}_\gamma$ on the closed geodesic $\gamma^*$ for each $\gamma
\in P_Y$.  We complete our model $N$ by adding a $Y$-cap: this is a
homotopy $$f_Y \colon S \times I \to C(X,Y)$$ from
$f_{|G|}\vert_{\bdry^+ N_{|G|}}$ to the convex core boundary $g_Y
\colon Y_h \to \bdry^+ C(X,Y).$   Gluing this homotopy $S\times I$,
to
$\bdry^+ N_{|G|}$ and extending $f_{|G|}$ over the $Y$-cap by $f_Y$,
we obtain the final piece of our model $N$ and the resulting map
$$f \colon N \to C(X,Y),$$ a homotopy equivalence whose
restrictions
$$f \vert_{\bdry^+ N} \to \bdry^+ C(X,Y)
\ \ \ \text{and} \ \ \  
f \vert_{\bdry^- N} \to \bdry^- C(X,Y)$$ are homeomorphisms.

Though the boundary $\bdry C(X,Y)$ is not generically smooth, by
taking the boundary
$\bdry \calN_\epsilon(C(X,Y))$ of the $\epsilon$-neighborhood of the
convex core we obtain a pair of $C^1$ surfaces $\bdry^+
\calN_\epsilon(C(X,Y))$ and $\bdry^- \calN_\epsilon(C(X,Y))$ with
$C^1$ path metrics
\cite[Lem. 1.3.6]{Epstein:Marden:convex}.  In the interest of
computing volume, we perturb $f$ to a piecewise smooth map
$$f^\epsilon \colon N \to \calN_\epsilon(C(X,Y))$$ 
by adjusting $f$ by a homotopy that changes $f$ only
on  ${\rm cap}_X$ and ${\rm cap}_Y$, so that 
$$
f^\epsilon \vert_{\bdry^+
N} \colon \bdry^+ N \to \bdry^+ \calN_\epsilon(C(X,Y))
\ \ \ \text{and} \ \ \ 
f^\epsilon \vert_{\bdry^- N} \colon \bdry^- N \to \bdry^-
\calN_\epsilon(C(X,Y)) $$ are homeomorphisms, and
$f^\epsilon$ is $C^1$ on the interiors of the caps of $N$.  The map $f$ is
already simplicial on the triangulated part of $N$ so $f^\epsilon$ is
piecewise $C^1$.  A degree argument shows  $f^\epsilon$ is  surjectve,
proving the theorem. 
\qed

\subsection{Bounding the volume}
Given a piecewise differentiable 3-chain $C$ in a hyperbolic
$3$-manifold $M$, 
the function $\deg_C \colon M \to \zed$ which measures the degree of $C$
in $M$ is well defined at almost
every point of $M$.  
We define the
{\em mass} $\mass(C)$ of $C$ to be the integral 
$$\mass(C) = \int_M |\deg_C| dV$$ where $dV$ is the hyperbolic
volume form on $M$ (cf. \cite[\S4]{Thurston:hype1}).

Moreover, if $F \colon P \to M$ is a map of a piecewise differentiable 
3-manifold $P$ to $M$, and $C$ is a piecewise differentiable 3-chain
in $P$, then we define the {\em $F$-mass} of $C$ by the integral
$$\mass_F(C) = \int_M |\deg_{F(C)}| dV.$$
The $F$-mass of $C$ bounds the volume
$\vol(F(C))$ of the image of $C$ in $M$.  
Hence, given our piecewise differentiable surjective map 
$$f^\epsilon \colon N \to 
\calN_\epsilon(C(X,Y))$$
the volume $\vol(\calN_\epsilon(C(X,Y))$, which bounds $\vol(X,Y)$, is
bounded by the sums of the  $f^\epsilon$-masses of the chains that
make up $N$.
In other words, if $N$ decomposes into 3-chains $C_k$, we have
$$\vol(X,Y) \le \sum_{C_k \subset N} \mass_{f^\epsilon}(C_k) =: \mass_{f^\epsilon}(N).$$

Thus, Theorem~\ref{theorem:volume:upper} will follow from the
following proposition.
\begin{prop}
Given $S$ there are constants $K_3 >1$ and
$K_4>0$ so that the map $f^\epsilon$ is properly homotopic to a map
$f_\theta^\epsilon
\colon N \to \calN_\epsilon(C(X,Y))$ for which 
$$\mass_{f^\epsilon_\theta}(N) \le K_3 d_{\bf P}(P_X,P_Y) + K_4.$$
\label{proposition:mass:upper}
\end{prop}

Let $\calV_3$ denote the maximal volume of a tetrahedron in
hyperbolic 3-space (see \cite[ch. 7]{Thurston:book:GTTM},
\cite{Benedetti:Petronio:book}). 
We begin our approach to Proposition~\ref{proposition:mass:upper}
with the following lemma. 

\begin{lem}There is a constant $K_\Delta$ so that 
the map $f^\epsilon$ is properly homotopic to a map
$f^\epsilon_\theta$ that also satisfies the conclusions of
Theorem~\ref{theorem:model} so that
$$\mass_{f^\epsilon_\theta}(N_\Delta) < K_\Delta \cdot \calV_3 \cdot
d_{\bf P}(P_X,P_Y) + 1.$$
\label{lemma:triangulated:bound}
\end{lem}

\bold{Proof:}  Because of the possibility of a large number of Dehn
twist blocks in $N_\Delta$, there is not in general a uniform constant
$K$ for which the number of tetrahedra used to triangulate $N_\Delta$
is less than $K d_{\bf P}(P_X,P_Y)$.  We will show, however, that by
modifying $f^\epsilon$ by a homotopy, we can force the tetrahedra that
lie in Dehn twist blocks to have mass as small as we like.  

By Theorem~\ref{theorem:model}, the
number of blocks in $N_\Delta$ that are not Dehn twist blocks is
bounded by $3 d_{\bf P}(P_X,P_Y)$.  Since there is a uniform bound to
the number of tetrahedra in any block, there is a constant $K_\Delta$
so that all but at most $K_\Delta d_{\bf P}(P_X,P_Y)$ 
tetrahedra of $N_\Delta$ lie in Dehn twist
blocks.

Let $\alpha^* \in {\rm spin}(G)$ be a spinning geodesic.  Let
$f^\epsilon_\theta$ 
be defined by the following homotopy of $f^\epsilon$ through maps that
are simplicial on $N_\Delta$: for each $\alpha \in \cup_j P_j$, slide
the vertices $p_\alpha$ and $\bar{p}_\alpha$ along the geodesic
$\alpha^*$ a distance $$\frac{\theta}{2\pi} \ell_{Q(X,Y)}(\alpha^*)$$
in the direction of the reference orientation chosen for $\alpha$.
(See \cite{Thurston:hype1} for another example of this {\em spinning}
of triangulations).  The following lemma shows that tetrahedra that
lie in a Dehn twist block, which we will call {\em Dehn twist
tetrahedra}, can be made to have small $f^\epsilon_\theta$-mass by
spinning to sufficiently high values of $\theta$.

\begin{lem}
If $\Delta \subset N_\Delta$ is a Dehn twist tetrahedron, then
$$\mass_{f^\epsilon_\theta}(\Delta) \to 0 \ \ \ \text{as} \ \ \ \theta 
\to \infty.$$
\label{lemma:tetraspin}
\end{lem}

\bold{Proof:}  Recall from Theorem~\ref{theorem:model}, each Dehn
twist tetrahedron $\Delta$ has at 
least one edge $e$ for which $f^\epsilon(e) \subset \alpha^*$ for some
spinning geodesic $\alpha^*$.

Lift $f_\theta^\epsilon$ to $\wt{f_\theta^\epsilon} \colon \wt{N} \to
\wt{\half^3}$ and choose a lift $\wt{\Delta}$ to $\wt{N}$.  Let
$\wt{\alpha^*}$ denote the lift of $\alpha^*$ to $\half^3$ for which
$\wt{f_\theta^\epsilon}$ sends the lifted edge $\wt{e} \subset
\wt{\Delta}$ to $\wt{\alpha^*}$.  Let $e'$ be the opposite edge of
$\Delta$  ($e$ and $e'$ have no endpoints in common).   

Let $I_\theta$ be the ideal tetrahedron in $\half^3$ for which
\begin{enumerate}
\item $\wt{f^\epsilon_\theta}(\wt{e'})$ lies in one edge $e_\infty'$
of $I_\theta$,
\item the two other edges $e_\infty^1$ and $e_\infty^2$  of $I_\theta$ emanating from one endpoint of
$e_\infty'$ pass through endpoints $p_1$ and $p_2$ of
$\wt{f^\epsilon_\theta}(\wt{e})$.
\end{enumerate}
The image $\wt{f^\epsilon_\theta}(\wt{\Delta})$ lies in $I_\theta$
(see Figure~\ref{figure:tetraspin}).
\makefig{Spinning
tetrahedra in
$\half^3$.}{figure:tetraspin}{\psfig{file=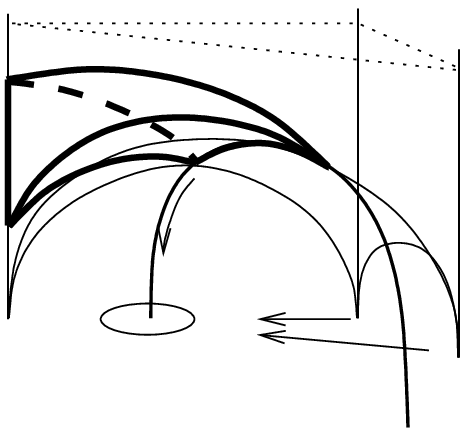,height=2.3in}}

Normalize by an isometry of $\half^3$ so that the edge $e_\infty'$ of
$I_\theta$ has its ideal endpoints at $0$ and $\infty$, and so that
the terminal fixed point of $\wt{\alpha^*}$ (with respect to the
reference orientation for $\alpha$) lies at $1 \in \cx$.

Then for any $r >0$ there exists a $\theta$ so that the two ideal
vertices of $I_\theta$ that do not lie at $0$ and $\infty$ lie within
a small disk $|1-z| < r$.  It follows that the dihedral angle of
$I_\theta$ at $e_\infty'$ tends to $0$ as $\theta \to \infty$.   
But the volume of an ideal tetrahedron tends to 0 as any of its
dihedral angles tends to 0, so we have
$$\mass_{f^\epsilon_\theta}(\Delta) < \vol(I_\theta) \to 0 \ \ \
\text{as} \ \ \ \theta \to \infty.$$
\qed

\ital{Continuation of the proof of Lemma~\ref{lemma:triangulated:bound}:} 
By Lemma~\ref{lemma:tetraspin}, for any Dehn twist block $B$, the
quantity  $\mass_{f^\epsilon_\theta}(B)$ is
as small as we like for $\theta$ sufficiently large.  
If $\calB^{\rm tw}$ denotes the union of all Dehn twist blocks in
$N_\Delta$, then, we may choose $\theta$ sufficiently large so that
$$\mass_{f^\epsilon_\theta}(\calB^{\rm tw}) < 1.$$

Since $f^\epsilon_\theta$ is simplicial on $N_\Delta$, we have
$$\mass_{f^\epsilon_\theta}(\Delta) < \calV_3$$ for any tetrahedron $\Delta
\subset N_\Delta$, and since all but at most $K_\Delta d_{\bf P}(P_X,P_Y)$
tetrahedra in $N_\Delta$ lie in $\calB^{\rm tw}$, we have
$$\mass_{f^\epsilon_\theta}(N_\Delta) < K_\Delta \cdot \calV_3 \cdot
d_{\bf P}(P_X,P_Y) + 1.$$
\qed

\bold{Bounding the volume of the caps.}
The bound on the $f^\epsilon_\theta$-mass of the triangulated part
$N_\Delta$ in terms of the distance $d_{\bf P}(P_X,P_Y)$ will be
sufficient for Proposition~\ref{proposition:mass:upper} once we show the following
uniform bound on the $f^\epsilon_\theta$-mass of the caps of $N$.

\begin{lem}
There is a uniform constant $K_{\rm cap}$, depending only on $S$ so
that after modifying $f^\epsilon_\theta$ by a homotopy on ${\rm
cap}_X$ we have 
$$\mass_{f^\epsilon_\theta}({\rm cap}_X) < K_{\rm cap}$$
and similarly for ${\rm cap}_Y$.
\label{lemma:cap:bound}
\end{lem}

\bold{Proof:}
Our goal will be to modify the homotopy $f^\epsilon_\theta \vert_{{\rm 
cap}_X}$ from $\bdry^-
\calN_\epsilon(C(X,Y))$ to the simplicial hyperbolic surface $Z_X$ 
by cutting the surface $S$ into annuli and controlling the trace of
the homotopy on each annulus (a solid torus).  To obtain such control, 
we choose this decomposition compatibly with the pants decomposition $P_X$.

We fix attention on a single pair of pants $\wh{S} \subset S- P_X$.
By a {\em figure-8 curve} on $\wh{S}$ we will mean a closed curve that
intersects itself once on $\wh{S}$ and divides $\wh{S}$ into three
annuli, one parallel to each boundary component of $\wh{S}$ (see
Figure~\ref{figure:figure8}).  To prove Lemma~\ref{lemma:cap:bound} we
establish the following basic lemma in hyperbolic surface geometry.
(We continue to treat the closed case; the case when $S$ has boundary
is similar).
\begin{lem}
Let $S$ a closed surface with negative Euler characteristic, and
let $L' \ge L$ be a constant greater than or equal to the Bers constant
$L$ for $S$.  Then there is a constant $L_8(L') >0$ so that the
following holds: if $P$ is a pants decomposition of $S$, $Z \in \Sing_k(S)$
is a possibly singular hyperbolic surface, and the length bound
$\ell_Z(\alpha) < L'$ holds for each $\alpha \in P$, then any
figure-$8$ curve $\gamma$ in any component $\wh{S} \subset S - P$
satisfies $$\ell_Z(\gamma) < L_8(L').$$
\label{lemma:figure8}
\end{lem}

\makefig{A figure eight curve.}{figure:figure8}{\psfig{file=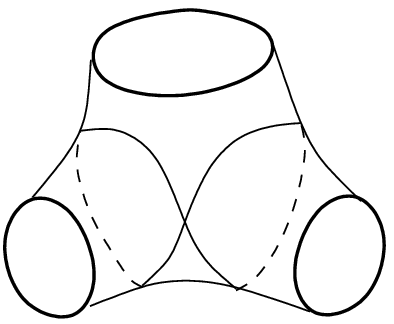,height=1.5in}}
\bold{Proof:}  Note that any bound on the geodesic length of a given
figure-8 curve $\gamma$ guarantees a bound on the geodesic length of
any other, by taking twice the original bound.

Assume first that $Z \in \Teich(S)$, so $Z$ is a non-singular
hyperbolic surface.
From the thick-thin decomposition for hyperbolic surfaces, given any
$\delta >0$, there is a uniform bound $D_\delta$ 
to the diameter $\diam(Z_{\ge \delta})$ of the $\delta$-thick part of
$Z$, where $D_\delta$
depends only on $\delta$ and the surface $S$.
Let $\bdry \wh{S} = \alpha_1 \disjunion \alpha_2
\disjunion \alpha_3$.  Since $\ell_Z(\alpha_i) < L'$, there is a
$\delta$ 
so that $\alpha^*_i$ can
only intersect $Z_{< 
\delta}$ in a component of $Z_{< \delta}$ for which it is 
the core curve.    Choosing $\delta$ smaller if necessary, we can
ensure that the boundary of any component of $Z_{< \delta}$ is a
curve of length less than $L'$.

Let $\wh Z$ be the realization of $\wh S$ as a subsurface of $Z$
bounded by the curves $\cup \alpha_i^*$. 
Either two of the $\alpha_i^*$ lie entirely in $Z_{\ge
\delta}$ or two of the $\alpha_i$ are homotopic into
$Z_{<\delta}$.  Without loss of generality, assume $\alpha_1^*$ and
$\alpha_2^*$ lie in $Z_{\ge
\delta}$.  Then they can be joined by an arc $b$ in $\wh Z$ of length
less than 
$D_\delta$.  Either $\alpha_1 \cdot b \cdot \alpha_2 \cdot b^{-1}$ or
$\alpha_1 \cdot b \cdot \alpha_2^{-1} \cdot b^{-1}$ is a figure-8
curve of length less than $2L' + 2D_\delta$.  Otherwise, if $A_1$ and
$A_2$ are components of $Z_{< \delta}$ representing the
homotopy classes of $\alpha_1$ and $\alpha_2$, then there
are components $a_1$ of $\bdry A_1$ and $a_2$ of $\bdry A_2$ with length
less than $L'$, and an arc $b$ joining $a_1$ to $a_2$  in
$\wh{Z}$ of length less than $D_\delta$.  By the same reasoning, there
is a figure-8 curve $\gamma$ of length $$\ell_Z(\gamma) < 2L' +
2D_\delta.$$

To treat the potentially singular case, let $Z \in \Sing_k(S)$ 
satisfy $\ell_Z(\alpha_i) < L'$ for each $\alpha_i \in P$.  Let
$Z^{\rm h} \in \Teich(S)$ represent the hyperbolic surface in the same
conformal class as $Z$.  If each $\alpha_i \in P$ satisfies
$$\ell_{Z^{\rm h}}(\alpha_i) < L'$$ then there is a figure-8 curve $\gamma
\subset \wh{S}$ with $\ell_{Z^{\rm h}}(\gamma) < 2L' + 2D_\delta$ by
the above reasoning.  By Ahlfors' lemma, \cite{Ahlfors:simplicial} we
have $$\ell_{Z}(\gamma) < \ell_{Z^{\rm h}}(\gamma) < 2L' + 2D_\delta$$
proving the lemma in this case.

If, however, some $\alpha_i \subset \bdry \wh{S}$ satisfies
$\ell_{Z^{\rm h}}(\alpha_i) \ge L' \ge L$, then by
Theorem~\ref{theorem:bers:constant} there is a simple closed curve
$\beta$, with $i(\alpha_i, \beta) \not= 0$, for which $\ell_{Z^{\rm
h}}(\beta) < L$.  Again applying Ahlfors' lemma, the curve $\beta$ has
length $$\ell_Z(\beta) <L$$ on $Z$, so its geodesic representative
$\beta^*$ furnishes an arc $b$ in $\wh{S}$ that joins $\alpha_i^*$ to
$\alpha_j^*$, or joins $\alpha_i^*$ to itself.  In either case, two
copies of $b$ together with arcs in the geodesics at its endpoints can
be assembled to form a figure-8 curve $\beta \subset
\wh{S}$ with length bounded by
$$\ell_Z(\beta) < 2L' + 2L \le 4L'.$$
By setting 
$$L_8(L')  = 4L' + 2D_\delta$$
the lemma follows.
\qed

\bold{Remark:}  Though we continue to work in the setting where $S$ is
closed, the proof of the Lemma for the case when $\bdry S \not=
\nullset$ goes through simply.

\ital{Continuation of the proof of Lemma~\ref{lemma:cap:bound}:} Again
consider the subsurface $\wh{S}$ and its realization as a 
subsurface $\wh{Z_X} \subset Z_X$ bounded by geodesics
$\hat{\alpha_1}$, $\hat{\alpha_2}$ and $\hat{\alpha_3}$ on $Z_X$.  
Let $$X_h^\epsilon = \bdry^+
\calN_\epsilon(C(X,Y))$$ denote the boundary component of
$\calN_\epsilon(C(X,Y))$ facing $X$, and let $\wh{X_h^\epsilon}$
denote the realization of $\wh{S}$ as a subsurface $X_h^\epsilon$
bounded by geodesics $\bar{\alpha_1}$, $\bar{\alpha_2}$ and
$\bar{\alpha_3}$ in the path metric on $X_h^\epsilon$.  Let $\gamma$
be a figure-8 curve on $\wh{S}$ with geodesic representatives
$\hat{\gamma}$ and $\bar{\gamma}$ on $\wh{Z_X}$ and
$\wh{X_h^\epsilon}$ respectively.

Let $A_1$, $A_2$ and $A_3$ denote the three annuli in $\wh{S} - \gamma$ 
for which $\alpha_i \subset \bdry A_i$, where $i = 1,2,3$, and let
$\gamma_i \subset \gamma$ be the loop on $\gamma$ so that $\gamma_i
\subset \bdry A_i$ for $i = 1, 2, 3$.  Let $\hat{A_i}$ be the
realization of $A_i$ on $\wh{Z_X}$ with $\bdry \hat{A_i} \subset
\hat{\alpha_i} \cup \hat{\gamma}$, and let
$\bar{A_i}$ be the realization of $A_i$ on $\wh{X_h^\epsilon}$ with
$\bdry \bar{A_i} \subset \bar{\alpha_i} \cup \bar{\gamma}$. 

By a theorem of D. Sullivan \cite{Sullivan:Bourbaki} elaborated upon by
D. Epstein and A. Marden \cite[Thm. 2.3.1]{Epstein:Marden:convex}
the path metric on the surface 
$X_h^\epsilon$
has bounded distortion from the hyperbolic
metric on $X$: the {\em nearest point retraction} map (see
\cite{Epstein:Marden:convex}) is $4 \cosh(\epsilon)$-Lipschitz, in
particular.  It follows that there is a constant $C_5 >4$ so that
choosing $\epsilon$ sufficiently small we have
$$\ell_{X_h^\epsilon}(\bar{\alpha_i}) < C_5 L$$
and 
$$\area(\wh{X_h^\epsilon}) < \area(X_h^\epsilon) < C_5 2 \pi |\chi(S)|$$
where $i = 1,2,3$.
Since $\wh{Z_X}$ is triangulated by $8$ hyperbolic triangles, we have
$$\area(\hat{A_i}) < 8 \pi$$
for $i = 1, 2, 3$.

Returning to the map $f^\epsilon_\theta \colon N \to
\calN_\epsilon(C(X,Y))$, we may define the $f^\epsilon_\theta$-mass of
a piecewise differentiable $2$-cycle in $N$ similarly to that of a
3-chain, by integrating the absolute value of the degree with respect
to two dimensional Lebesgue-measure on $M$
(cf. \cite[Prop. 4.1]{Thurston:hype1} and the preceding discussion).

We modify the map $f_\theta^\epsilon$ on ${\rm cap}_X \cong S \times
I$ as follows.  By Theorem~\ref{theorem:model}, the map
$f_\theta^\epsilon$ already sends $\alpha_i 
\times 0$ to the geodesic $\hat{\alpha_i}$ for each $\alpha_i \in
P_X$.  Straighten $f_\theta^\epsilon$ without changing its values on
$\bdry^+ {\rm cap}_X$ so that $f^\epsilon_\theta$ sends each $\alpha_i
\times \{0\}$ to the geodesic $\bar{\alpha_i} \subset X_h^\epsilon$,
and so that $f^\epsilon_\theta$ sends each track $x \times [0,1]$, 
where $x$ lies in $\alpha_i$ to a geodesic.  The image
$f^\epsilon_\theta (\alpha_i \times [0,1])$ is a ruled annulus, and
any ruled annulus has area less than the sum of the lengths of its
boundary components (see e.g. \cite[Ch. 9]{Thurston:book:GTTM}
\cite[\S 3.2]{Bonahon:tame}).
Since 
$\ell_{X_h^\epsilon}(\bar{\alpha_i}) < C_5L$, and since $\hat{\alpha_i}$ 
is the geodesic representative of $\alpha_i$ in $Q(X,Y)$, we have
$$\mass_{f^\epsilon_\theta} (\alpha_i \times [0,1]) < 2C_5L.$$

Applying Lemma~\ref{lemma:figure8} to $X$ with $L'=L$, 
we have an $L_8 = L_8(L)$
so that
$$\ell_X(\gamma) < L_8.$$
It follows that
$$\ell_{X_h^\epsilon}(\gamma) < C_5 L_8.$$ Applying
Lemma~\ref{lemma:figure8} to $Z_X$ with $L' = C_5L$, we have a
$L_8' = L_8(C_5 L)$ for which 
$$\ell_{Z_X}(\gamma) < L_8'$$
Letting $$K_8 = \max\{L_8, L_8'\},$$
we may straighten $f^\epsilon_\theta$ to send $\gamma \times \{0\}$ to
its geodesic representative $\bar{\gamma} \subset X_h^\epsilon$ in
$X_h^\epsilon$ and straighten further as before so that for $x \in
\gamma$ the map $f^\epsilon_\theta$ sends $x \times [0,1]$ to a
geodesic.  Arguing as for $\alpha_i$, each annulus $\gamma_i \times
[0,1]$ has mass $$\mass_{f^\epsilon_\theta}(\gamma_i \times [0,1]) < 2
K_8.$$

It follows that the union of the four annuli 
$$\eT_i = \left(A_i \times \{1\}\right) \cup \left(\alpha_i \times
[0,1]\right) \cup \left( A_i \times \{0\} \right)
\cup \left( \gamma_i \times [0,1] \right)$$
along their boundaries is a torus with mass 
$$\mass_{f^\epsilon_\theta}(\eT_i) < 2 K_8 + (8 + 2|\chi(S)|) \pi = K_8'.$$

Applying \cite[Prop. 4.1]{Thurston:hype1}, 
we have that the solid torus $\eV_i \subset N$ bounded by $\eT_i$ has mass
$$\mass_{f^\epsilon_\theta}(\eV_i)  \le
\mass_{f^\epsilon_\theta}(\eT_i)< K_8'.$$ 
It follows that the total $f^\epsilon_\theta$-mass of $\wh{S}\times I$
is less than $3 K_8'$, and since the number of pieces of the complement 
$S - P_X$ depends only on $S$, it follows that 
$$\mass_{f^\epsilon_\theta}({\rm cap}_X) < K_{\rm cap}$$
for an a priori constant $K_{\rm cap}$.

\qed

\bold{Conclusion.}
The proof of Proposition~\ref{proposition:mass:upper} 
is now an application of the preceding lemmas.

\bold{Proof:} {\em (of Proposition~\ref{proposition:mass:upper}).}  Applying
lemmas~\ref{lemma:triangulated:bound} and~\ref{lemma:cap:bound}, we
have
\begin{eqnarray*}
\mass_{f^\epsilon_\theta}(N) 
&=& \mass_{f^\epsilon_\theta}(N_\Delta) + \mass_{f^\epsilon_\theta}({\rm
cap}_X) + \mass_{f^\epsilon_\theta}({\rm cap}_Y) \\
&<& K_\Delta \cdot \calV_3 \cdot d_{\bf P}(P_X,P_Y) + 1 + 2K_{\rm cap}.
\end{eqnarray*}
Setting $K_\Delta \cdot \calV_3 = K_3$ and $1+ 2 K_{\rm cap} = K_4$,
the result follows.

\ital{The case when $\bdry S \not= \nullset$.}  Minor modifications
to the above arguments are required when $\bdry S \not= \nullset$, and
$Q(X,Y)$ has peripheral rank-1 cusps.
The primary difference in our model $N$ in this case arises from the fact
that boundary curves of $S$ are not elements of the pants
decompositions $P_j$ interpolating between $P_X$ and $P_Y$.  For the
purposes of the construction, triangulations suited to $P_j$ are then
triangulations of $S$ with two vertices on each boundary component so
that on each pair of pants in $S - P_j$ the triangulation has the
structure of Definition~\ref{defn:triangulation}.  The blocks required
to accomplish each elementary move are then identical; to map them in
to $C(X,Y)$ we need only specify the behavior of the maps $f_{j,*}$ on
the boundary curves.

Indeed, if $\gamma \subset \bdry S$, there is no
geodesic in the free homotopy class of $\gamma$.
Choosing a horocycle $h$ in
the homotopy class of $\gamma$ on the boundary $\bdry (C(X,Y)_{\ge
\epsilon_0})$ and a
pair of antipodal vertices $p_h$ and $\bar{p}_h$ on $h$, we let
$\gamma^*$ denote the piecewise geodesic obtained from straightening
$h$ {\em rel}-$p_h \cup
\bar{p}_h$.  The behavior of $f_{j,*}$ is then determined by
requiring that $f_{j,*}$ send each $\gamma \subset \bdry S$ to the
corresponding $\gamma^*$ and that $f_{j,*}$ be simplicial as before.

If $f_\Delta \colon N_\Delta \to C(X,Y)$ denotes the simplicial map of
$N_\Delta$ into $C(X,Y)$ given by the above instructions, we modify
$f_\Delta$ by sending $\epsilon_0$ to $0$: the edges and vertices of
$N_\Delta$ incident on $\gamma \subset \bdry S$ are mapped by
$f_\Delta$ deeper and deeper into the $\gamma$-cusp so that the
limiting map $f_{\Delta,0}$ sends each such tetrahedron to a
tetrahedron or triangle with one vertex at infinity.  The restriction
of $f_{\Delta,0}$ to the interiors of $\bdry^+ N_\Delta$ and $\bdry^-
N_\Delta$ factor through simplicial hyperbolic surfaces $h_X \colon
Z_X \to C(X,Y)$ and $h_Y \colon Z_Y \to C(X,Y)$ realizing $P_X$ and
$P_Y$, and the restriction of $f_{\Delta,0}$ to $N_\Delta$ with the
$\gamma$ curves removed is a proper homotopy from $h_X$ to $h_Y$.

The caps, now proper homotopies from $Z_X$ to $X_h^\epsilon$ and
$Z_Y$ to $Y_h^\epsilon$ may then be added to the triangulated part.
The total mass is again bounded by $K_{\rm cap}$, since the caps still
decompose into solid tori with boundary of bounded area (some are now
bounded by properly embedded annuli asymptotic to cusps).  The
resulting map, a proper homotopy between the two components
$X_h^\epsilon$ to $Y_h^\epsilon$ of $\bdry \calN_\epsilon(C(X,Y))$ may
be spun about the spinning geodesics as before to force the Dehn twist
tetrahedra to have arbitrarily small mass.
Spinning sufficiently far, then, we again have the conclusion of
Proposition~\ref{proposition:mass:upper}.

\qed

\section{Geometric limits}

We conclude with the following application of our results to the
study of algebraic and geometric limits of hyperbolic 3-manifolds.

\begin{theorem}
Let $\{Q(X_k,Y_k)\}_{k=1}^\infty \subset AH(S)$ be an algebraically
and geometrically convergent 
sequence with geometric limit $N_G$.  Then $N_G$ is geometrically
finite if and only if there is a $K >0$ for which 
$$d_{\rm WP}(X_k,Y_k) <K$$ for all $k$.
\end{theorem}

\bold{Proof:}  If the geometric limit $N_G$ is geometrically infinite, 
the volumes of the convex cores of the approximates $Q(X_k,Y_k)$ grow without 
bound (see \cite[Lem. 7.1]{Canary:Minsky:tame:limits}), which,
applying Theorem~\ref{theorem:main}, implies
that $d_{\rm WP}(X_k,Y_k)$ grows without bound.

Likewise, applying \cite{Taylor:thesis} or
\cite[Thm. 3.1]{McMullen:HDI}, when $N_G$ is geometrically finite the
$\epsilon$-thick parts of the convex cores of $Q(X_k,Y_k)$ converge
geometrically to that of $N_G$ for $\epsilon$ sufficiently small.  It
follows that the volume of $\core(Q(X_k,Y_k))_{\ge \epsilon}$
converges to the volume of $\core(N_G)_{\ge \epsilon}$ for each
sufficiently small $\epsilon$ and thus that $$\vol(\core(Q(X_k,Y_k)))
\to \vol(\core(N_G))$$ as $k$ tends to $\infty$.  The theorem follows.
\qed

\bibliographystyle{math}
\bibliography{math}

\end{document}